\documentclass[12pt]{amsart}
\usepackage{hyperref}
\begin{document}
\numberwithin{equation}{section}

\def\1#1{\overline{#1}}
\def\2#1{\widetilde{#1}}
\def\3#1{\widehat{#1}}
\def\4#1{\mathbb{#1}}
\def\5#1{\frak{#1}}
\def\6#1{{\mathcal{#1}}}

\def\C{{\4C}}
\def\R{{\4R}}
\def\N{{\4N}}
\def\Z{{\4Z}}
\def\Q{{\4Q}}

\title[]{A geometric approach to Catlin's boundary systems}
\author[D. Zaitsev]{Dmitri Zaitsev}
\address{D. Zaitsev: School of Mathematics, Trinity College Dublin, Dublin 2, Ireland}
\email{zaitsev@maths.tcd.ie}

\begin{abstract}
For a point $p$ in a smooth real hypersurface 
$M\subset\C^n$, 
where the Levi form has the nontrivial kernel $K^{10}_p$,
we introduce an invariant cubic tensor
$$\tau^3_p \colon \C T_p \times  K^{10}_p \times \1{K^{10}_p} 
\to \C\otimes (T_p/H_p),$$ 
which together with Ebenfelt's tensor $\psi_3$,
constitutes the full set of the $3$rd order invariants of $M$ at $p$.

Next, in addition, assume $M\subset\C^n$ to be 
 {\em (weakly) pseudoconvex}.
Then $\tau^3_p$ must identically vanish. 
In this case we further define an invariant
quartic tensor
$$\tau^4_p \colon  \C T_p \times \C T_p 
	\times K^{10}_p\times \1{K^{10}_p} \to \C\otimes (T_p/H_p),$$
	and for every $q=0, \ldots, n-1$, 
	an invariant {\em submodule sheaf}
	 $\6S^{10}(q)$ of $(1,0)$ vector fields
in terms of the Levi form,
and an invariant {\em ideal sheaf} $\6I(q)$ of complex functions
generated by certain components and derivatives of the Levi form,
such that the set of points of Levi rank $q$
is locally contained in real submanifolds
defined by real parts of the functions in $\6I(q)$,
whose tangent spaces have explicit algebraic description
in terms of the quartic tensor $\tau^4$.

Finally, we relate the introduced invariants
with D'Angelo's finite type, Catlin's mutlitype 
and Catlin's boundary systems.
\end{abstract}

\maketitle
\tableofcontents

\def\Label#1{\label{#1}}


\def\cn{{\C^n}}
\def\cnn{{\C^{n'}}}
\def\ocn{\2{\C^n}}
\def\ocnn{\2{\C^{n'}}}


\def\dist{{\rm dist}}
\def\const{{\rm const}}
\def\rk{{\rm rank\,}}
\def\id{{\sf id}}
\def\aut{{\sf aut}}
\def\Aut{{\sf Aut}}
\def\CR{{\rm CR}}
\def\GL{{\sf GL}}
\def\Re{{\sf Re}\,}
\def\Im{{\sf Im}\,}
\def\span{\text{\rm span}}
\def\ord{\text{\rm ord}}

\def\codim{{\rm codim}}
\def\crd{\dim_{{\rm CR}}}
\def\crc{{\rm codim_{CR}}}

\def\phi{\varphi}
\def\e{\varepsilon}
\def\eps{\varepsilon}
\def\d{\partial}
\def\a{\alpha}
\def\b{\beta}
\def\g{\gamma}
\def\G{\Gamma}
\def\D{\Delta}
\def\Om{\Omega}
\def\k{\kappa}
\def\l{\lambda}
\def\L{\Lambda}
\def\z{{\bar z}}
\def\w{{\bar w}}
\def\Z{{\1Z}}
\def\t{\tau}
\def\th{\theta}
\def\p{\phi}
\def\de{\delta}
\def\la{\langle}
\def\ra{\rangle}
\def\r{\rho}

\emergencystretch15pt
\frenchspacing

\newtheorem{Thm}{Theorem}[section]
\newtheorem{Cor}[Thm]{Corollary}
\newtheorem{Pro}[Thm]{Proposition}
\newtheorem{Lem}[Thm]{Lemma}

\theoremstyle{definition}\newtheorem{Def}[Thm]{Definition}

\theoremstyle{remark}
\newtheorem{Rem}[Thm]{Remark}
\newtheorem{Exa}[Thm]{Example}
\newtheorem{Exs}[Thm]{Examples}

\def\bl{\begin{Lem}}
\def\el{\end{Lem}}
\def\bp{\begin{Pro}}
\def\ep{\end{Pro}}
\def\bt{\begin{Thm}}
\def\et{\end{Thm}}
\def\bc{\begin{Cor}}
\def\ec{\end{Cor}}
\def\bd{\begin{Def}}
\def\ed{\end{Def}}
\def\br{\begin{Rem}}
\def\er{\end{Rem}}
\def\be{\begin{Exa}}
\def\ee{\end{Exa}}
\def\bpf{\begin{proof}}
\def\epf{\end{proof}}
\def\beq{\begin{equation}}
\def\eeq{\end{equation}}

\def\ben{\begin{enumerate}}
\def\een{\end{enumerate}}

\section{Introduction}
\subsection{Overview for broader audience}
In this brief overview we put this paper's material
in somewhat broader context.
The methods and tools introduced here
may be of interest for general systems of Partial Differential Equations (PDE),
beyong the context of the $\bar\d$-equation considered here.
An evidence for this is
the recent breakthrough paper by Siu~\cite{S17}
extending to general PDE systems 
the celebrated multiplier ideal technique 
of Kohn~\cite{K79}.

In the study of PDE systems and their solutions,
an important general approach is that of the {\em a priori estimates}
$\|u\|_1\le \|u\|_2$, 
where $\|\cdot\|_1, \|\cdot\|_2$ is a pair of semi-norms and
 $u$ is a test function.
Finding a priori estimates is typically a difficult problem
whose solution in most cases relies on the specific nature
of the sytem, with very few general approaches known.
The multiplier ideal technique by Kohn~\cite{K79}
and the potential-theoretic approach 
by Catlin~\cite{C84ann, C87},
are ones of the few known general approaches
for the $\bar\d$-Neumann problem,
the boundary value problem for the $\bar\d$-equation 
with Neumann boundary conditions,
see recent expositions in \cite{CD10, S10, St10, MV15, S17}.

In both approaches,  understanding singularities 
of the Cauchy-Riemann structure
induced on the boundary by the ambient complex structure,
is of utmost importance.
Our goal here is to develop new geometric invariants
to tackle this problem.
A particular advantage of geometric invariants
comes from the freedom of using them in arbitrary coordinates,
as well as providing certain adapted coordinates,
where the computations can be significantly simplified.

The invariants introduced are {\em tensors, ideal sheaves of functions
and submodule sheaves of vector fields}.
The tensors arise from derivatives
 of the Levi form,
a fundamental invariant of the induced Cauchy-Riemann structure.
However, only derivatives along certain special vector fields
lead to invariant tensors, as illustrated by examples in this paper.
This observation leads to the study of
suitable submodule sheaves,
where these special vector fields must be contained,
in order to define tensors.

While being a powerful tool for computations in arbitrary coordinates,
higher order tensors, as opposed to the Levi form, have the 
fundamental limitation of not well-behaving across singularities,
since they are defined on kernels of varying dimension.
In order to achieve some better behaviour and control,
more flexible objects are needed,
such as ideal sheaves of functions.
In this paper we introduce
invariant ideal sheaves
generalizing functions
in Catlin's boundary systems \cite{C84ann},
whereas the invariant tensors
control the transversality and 
nondegeneracy property
of those functions.

\subsection{Detailed overview}
%
%
%
%
%
%
In more concrete terms,
the goal of this paper
is to introduce new geometric invariants
giving insight into some techniques developed by Catlin
in his celebrated papers 
\cite{C84ann, C87},
following previous foundational work by Kohn \cite{K64a, K64b, K72, K79}
on the $\bar\d$-Neumann problem.
In particular, we introduce new invariant ideal sheaves
containing functions arising in Catlin's boundary systems,
and  new invariant tensors permitting
to simplify the interative construction of the boundary systems
by more direct computations of the tensors' kernels.
The obtained geometric approach may lead
to sharper subelliptic etimates,
as well as to advance our understanding of
the Kohn's multiplier ideal sheaves \cite{K79}, 
in view of their relation with Catlin's technique,
as indicated by Siu \cite{S05, S10, S17}
(see also Nicoara~\cite{N14}).
Furthermore,
in his recent fundamental work \cite{S17},
Siu proposed new techniques for generating multipliers 
for general systems of partial differential equations,
also including a new procedures even
for the case of the $\bar\d$-Neumann problem.

\subsection{Conditions of property (P) type}
The importance of having a better understanding is further 
underlined
by the role played by the Catlin's potential-theoretic ``Property (P)'' type conditions
(see e.g.\ \cite{FS01, St06, St10, MV15, BS16} for recent surveys)
that found vast applications in multiple research areas
such as:

\ben

\item Compactness of the Kohn's $\bar\d$-Neumann solution operator
by Henkin and Iordan~\cite{HI97},
McNeal~\cite{M02},
Raich and Straube~\cite{RS08},
Harrington~\cite{Ha07, Ha11},
\c Celik and \c Sahuto\u glu~\cite{CS12}.
It  was even proved to be equivalent to Property (P)
for Hartogs domains in $\C^2$ by Christ and Fu~\cite{CF05},

\item Subelliptic estimates by 
Forn\ae ss and Sibony~\cite{FS89},
Straube~\cite{St97},  
Herbig~\cite{H07},
Harrington~\cite{Ha07}.

\item Invariant metric estimates 
due to Catlin~\cite{C89}, Cho~\cite{Ch92, Ch94, Ch02},
Boas-Straube-Yu~\cite{BSY95}, McNeal~\cite{M01}
and Herbort~\cite{He14}
and via subelliptic estimates in \cite{M92a}.

\item Stein neighborhood bases by
Harrington~\cite{Ha08}
and \c Sahuto\u glu~\cite{Sa12}.

\item Estimates and comparison of the Bergman and Szeg\"o kernels by
Boas~\cite{B87},
Nagel, Rosay, Stein and Wainger and \cite{NRSW89},
Boas-Straube-Yu~\cite{BSY95},
Charpentier and Dupain~\cite{CD06a, CD06b, CD14}
 Chen and Fu~\cite{CF11},
 Khanh and Raich~\cite{KhR14}.
 
\item Holomorphic Morse inequalities
and eigenvalue asymptotics for $\bar\d$-Neumann Laplacian
by Fu and Jacobowitz~\cite{FJ10}.

\item Tangential $\bar\d_b$ and
complex Green operator
by Raich-Straube~\cite{RS08},
Raich~\cite{R10},
Straube~\cite{St12},
Khanh, Pinton and Zampieri~\cite{KPZ12}.

\item
Construction of peak and 
bumping functions by
Diederich and Herbort~\cite{DH94},
Forn\ae ss and McNeal~\cite{FM94},
by Yu~\cite{Y94},
Bharali and Stens\o nes~\cite{BhS09},
as well as some generalisations of Property (P) by Khanh-Zampieri~\cite{KhZa12}.

\item
Division problems for holomorphic functions with $C^\infty$ boundary values
by Bierstone and Milman~\cite{BM87},
as an application of global regularity,
whose proof for smooth finite type boundaries
relies on Catlin's method.

\item
Regularity of solutions to the complex Monge-Amp\`ere equation
by Ha and Khanh~\cite{HK15} and Baracco, Khanh, Pinton and Zampieri~\cite{BKPZ16}.

\een


For bounded pseudoconvex domains 
with {\em real-analytic} boundaries of finite type in $\C^n$,
Kohn's \cite{K79} celebrated theory of subelliptic multipliers
provides an alternative approach to Catlin's in establishing subelliptic estimates.
The same approach also yields subelliptic estimates
for smoothly bounded domains of finite type in $\C^2$,
that were already treated by Kohn in his earlier paper \cite{K72}.
However, for general {\em smoothly bounded domains of finite type},
it remains open at the time of writing, whether the multiplier approach yields subelliptic estimates,
with Catlin's method being currently the only available.


\subsection{Submanifolds containing multitype level sets}
A key geometric aspect of Catlin's subelliptic estimates proof
consists of showing the existence of the so-called {\em weight functions}
 satisfying certain boundedness and positivity estimates for their complex Hessians,
 that are known as ``Property (P)'' type conditions 
 (see e.g.\ \cite{MV15, BS16} for recent surveys).
A major difficulty when constructing such weight functions
under geometric conditions (such as finite type), 
is to keep the uniform 
nature of the estimates across points 
of {\em varying ``degeneracy''}
for the underlying geometry.
A simple example of a degeneracy measure
is the {\em rank of the Levi form} of the 
boundary $M:=\d D$ (where $D$ is a domain in $\cn$).
A more refined measure
is the Catlin multitype \cite{C84ann}, see also \S\ref{m-type}.
To deal with points of varying multitype, 
Catlin developed his machinery of boundary systems \cite{C84ann}.
The main idea
to gain control of the multitype level sets is
by including them locally into certain {\em ``containing submanifolds''}.
A result of this type is the content of 
\cite[Main Theorem, Part (2)]{C84ann},
where a containing submanifold
is constructed by a collection of 
inductively chosen {\em boundary system functions}
that arise as certain carefully selected
(vector field) derivatives of the Levi form.

In this paper we focus on geometric invariants
behind the containing submanifold construction,
with the goal to extend and simplify the boundary system approach.
Our main discovery is that at the $4$th order level,
the boundary systems, as well as the type and the multitype,
can be described in terms of the new invariant objects,
such as tensors, submodule and ideal sheaves.

At the $4$th order level, the multitype level sets
boil down to simpler {\em level sets of the
Levi (form) rank} (see Proposition~\ref{multi-quartic} for details).
Recall that Catlin's boundary system functions
\cite{C84ann}
are constructed inductively with every new equation
depending on chosen solutions for previous ones.
In comparison, we here collect  defining functions 
for the Levi rank level sets
into {\em invariant ideal sheaves $\6I(q)$} on $M$,
for each Levi rank $q$.
The sheaf $\6I(q)$ is generated
by certain $1$st order Levi form derivatives as
described in Theorem~\ref{main} below.
In particular, arbitrary defining functions from $\6I(q)$
can be combined without any additional relations.
Furthermore, additional derivatives
of the Levi form along arbitrary complex vector fields $L^3$
(in Theorem~\ref{main}, Part (5)),
including transversal ones,
 are allowed for functions in $\6I(q)$.
In comparison, for a related boundary system function
given by the same formula,
the outside vector field $L^3$
would have to be in a special subbundle
inside the holomorphic tangent bundle.
As a result, we obtain richer classes of defining functions
allowing for more control over containing submanifolds
(see Example~\ref{transv-f}),
that may potentially lead to sharper a priori estimates.

%

%


%

In parallel to the ideal sheaf $\6I(q)$ construction,
we introduce {\em invariant quartic tensors} $\t^4$,
giving a precise control over differentials of
the functions in $\6I(q)$.
This is expressed in Theorem~\ref{main}, part (2),
where the tangent space of the containing manifold $S$
equals the real kernel of the tensor.
Importantly, the full tangent space of $S$
(rather than only the tangential part) is controlled here
via the kernel of $\t^4$, which means that transversal vector fields
must also be allowed among tensor arguments.
The tensors are constructed in Lemma~\ref{d2}
as certain $2$nd order Levi form derivatives
taken along all possible vector fields.
In comparison, only derivatives with 
respect to $(1,0)$ and $(0,1)$ vector fields
can appear in the boundary systems.
%
%
%
%

For reader's convenience,
we summarize the main results and constructions
in Theorem~\ref{main},
leaving more detailed and general statements
with their proofs in the chapters following.

\subsection{More details on invariant tensors and ideal sheaves}
Our first step in defining invariant tensors is a byproduct result 
giving a {\em complete set of cubic invariants}
for a general real hypersurface $M$, without pseudoconvexity assumption.
This is achieved by constructing an invariant cubic tensor $\t^3$
obtained by differentiating the Levi form along vector fields
with values in the Levi kernel, see Lemma~\ref{levi-der}.
Remarkably, to obtain tensoriality of the Levi form derivatives, 
it is of crucial importance to require
{\em both vector fields inside the Levi form to take values in the Levi kernel}
as explained in Example~\ref{one-ker}.
This stands further, in remarkable contrast with the cubic tensor $\psi_3$
defined by Ebenfelt \cite{E-jdg}
(by means of the Lie derivatives of the contact form),
where only one of the arguments needs to be in the Levi kernel.
%
On the other hand, Ebenfelt's tensor $\psi_3$
does not allow for transversal directions as $\t^3$ does.
It turns out that the pair $(\psi_3, \t^3)$ yields
a complete set of cubic invariants, 
as demonstrated by a normal form 
(of order $3$)
in Proposition~\ref{3-normal} eliminating all other terms
that are not part of either of the tensors.

We also investigate the construction based on double Lie brackets
(also considered by Webster \cite{W95}).
This approach, however, in order to yield a tensor, has to require all vector fields
to be in the complexified holomorphic tangent bundle,
leading only to a restriction of the cubic tensor $\t^3$.
Again, the double Lie bracket construction
is only tensorial when both vector fields inside the inner bracket
take their values in the Levi kernel (see Example~\ref{one-ker}).

As mentioned earlier, the cubic tensor $\t^3$ is constructed without any pseudoconvexity assumption.
On the other hand, in presence of pseudoconvexity,
the whole tensor $\t^3$
must vanish identically (Lemma~\ref{psc-vanish}).
The only cubic terms that may survive are of the form \eqref{psc-cubic}
which can never appear in the lowest weight terms,
and hence never play a role in Catlin's multitype 
and boundary system theory.  

Motivated by the above, we next look for quartic tensors.
It turns out (Example~\ref{ex-4}) 
that this time, neither second order Levi form derivatives
nor quartic Lie brackets provide tensorial invariants
even when all vector field arguments take their values in the Levi kernel.
To overcome this problem, we {\em restrict the choice of 
the vector fields involved} by requiring a certain kind of condition of
 ``Levi kernel inclusion up to higher order'' (Definition~\ref{ker-1}).
 In Lemma~\ref{1-ker-def} we show that
this additional condition always holds
for any vector field that is Levi-orthogonal
to a maximal Levi-nondegenerate subbundle,
which, in particular, arises in Catlin's boundary system construction.
However, the mentioned Levi-orthogonality lacks some invariance
as it depends on the choice of the subbundle.
In contrast, the Levi kernel inclusion up to order $1$
is invariant and only depends
on the $1$-jet of the vector field at the reference point.

With that restriction in place, an invariant quartic tensor $\t^4$
can now be defined in a similar fashion.
Then
it's restriction $\t^{40}$ enters 
the lowest weight normal form 
with weights $\ge 1/4$,
see Proposition~\ref{4-normal}.
It turns out, the restriction $\t^{40}$ 
provides exactly the missing information 
at the lowest weight level for hypersurfaces of finite type $4$
(where finite type $3$ cannot occur for pseudoconvex points in view of
Corollary~\ref{pse-nf-cor}).
For example, both D'Angelo finite type $4$
and Catlin's multitype up to entry $4$
can be completely characterized in terms
of $\t^{40}$.
In fact, having the finite type $4$ is equivalent
to the nonvanishing of $\t^{40}$ on complex lines
(Proposition~\ref{type-quartic}),
whereas having a multitype up to entry $4$
is equivalent to $\t^{40}$ having trivial kernel
(Propositions~\ref{multi-quartic}).

In \S\ref{ideal-par} we use the quartic tensor $\t^4$ to 
characterize the differentials of functions in 
the ideal sheaf $\6I(q)$
as well as the minimal tangent spaces
of containing manifolds
defined by a transversal set of functions in $\6I(q)$.

Finally, in \S\ref{bs} we obtain a characterization 
for a Catlin's boundary system,
where the most difficult part 
of obtaining vector field directions of nonvanishing Levi form derivatives
at the lowest weight 
is replaced by the nonvanishing of the tensor $\t^{40}$
on the vector fields' values at the reference point,
a purely algebraic condition.

In a forthcoming paper will shall 
extend the present 
geometric approach 
towards its {\em approximate versions}
with necessary control
to perform the induction step
in the subelliptic estimate proof.

\bigskip

{\bf Acknowledgements.}
The author would like to thank
J.J. Kohn, D.W. Catlin, J.P. D'Angelo,
E.J. Straube, M. Kolar, S. Fu,
J.D. McNeal and A.C. Nicoara
for numerous inspiring discussions.

\section{Notation and main results}
\Label{results}
We shall work in the smooth ($C^\infty$) category
unless stated otherwise.
Let $M\subset \C^{n}$ be a (smooth) real hypersurface.
We write $T:= TM$ for its tangent bundle, 
$H = HM \subset  T$ for the complex (or holomorphic) tangent bundle,
$Q:= T/H$ for the normal bundle,
as well as 
$$
	\C T := \C\otimes T,
	\quad 
	\C H: =\C\otimes H,
	\quad
	\C Q:= \C\otimes Q,
$$
for their respective complexifications.
Further, $H^{10}$ and $H^{01}= \1{H^{10}}$ denote $(1,0)$ and $(0,1)$ bundles respectively, such that $\C H = H^{10}\oplus \1{H^{10}}$.
By a small abuse of notation, we write $L\in V$
whenever a vector field $L$ is a local section in $V$,
where $V$ can be a bundle or a sheaf.
Similarly we write $V_1\subset V_2$
when $V_1$ is a local subbundle or subsheaf of $V_2$,
where it will be convenient to allow 
$V_1$ to be defined over smaller open sets
than those where $V_2$ is defined.

A subscript $p$ always denotes
evaluation at a point $p\in M$,
i.e.\ $L_p$ for the value of a vector field $L$ at $p$,
or $V_p$ for the space of all values
of elements in $V$,
which is the fiber when $V$ is a vector bundle.

On the dual side, 
$\Om = \Om M$ stands for the bundle of all real $1$-forms on $M$, 
$\C \Om$ for all complex $1$-forms,
$\Om^0$ for all {\em contact forms}, 
i.e.\ forms $\Om$ that are vanishing on $H$ and real-valued on $T$,
and $\C\Om^0$ for the corresponding complexification.

Recall that a {\em (local) defining function} of $M$
is any real-valued function $\rho$ in a neighborhood in $\C^n$
of a point in $M$,
with $d\rho\ne 0$ such that $M$ given by 
$\rho=0$ in the domain of $\rho$. 
For any defining function $\rho$, the one-form 
$\th: = i\d\rho$ 
spans (over $\R$) the bundle $\Om^0$ of all contact forms. 
 
We shall consider the standard 
$\C$-bilinear
pairing $\la \th, L\ra := \th(L)$
for $\th\in \C\Om$, $L\in \C T$.
By a slight abuse, we keep the same notation also for the induced pairing
$$
	\la \cdot, \cdot \ra \colon \C\Om^0 \times \C Q \to \C
$$
between the (complex) contact forms and the normal bundle.
With this notation, we regard the {\em Levi form tensor} at a point $p\in M$ 
as the  $\C$-bilinear map
$$
	\t^2_p \colon H^{10}_p \times \1{H^{10}_p} \to \C Q_p,
$$
which
is uniquely determined by the identity
\beq\Label{levi-id}
	\la \th_p, \t^2_p(L^2_p, L^1_p) \ra = i \la \th, [L^2, L^1] \ra_p,
	\quad 
	L^2\in H^{10}, \, L^1\in \1{H^{10}},
\eeq
where, as mentioned before, 
the membership notation for $L^2, L^1$ (such as $L^2\in H^{10}$)
means being local sections of the
corresponding bundles.
The normalization of $\t^2$ used here is chosen
such that for the quadric
$$
	\rho = -2\Re w + q(z,\z) =0,
	\quad
	(w,z)\in \C\times \C^{n-1},
$$ 
with $(1,0)$ vector fields
$$
	L_j := \d_{z_j} + q_{z_j }\d_w,
	\quad j=1,2,
$$
and the contact form 
$\th = i\d \rho = i(-dw + \d q) $,
we have
$$
	\la \th_0,  \t^2_0(\d_{z_j}, \d_{\z_k}) \ra 
	= i \la 
		\th, 
		[L_j, \1{L_k}]
	\ra_0
	= - \la
		dw,
		q_{z_j\z_k}(\d_{\w}-\d_w)
	\ra_0
	= \d_{z_j} \d_{\z_k} q,
$$
or more generally
\beq\Label{levi-calc}
	\la \th_0,  \t^2_0(v^2, v^1) \ra = \d_{v^2} \d_{v^1} q,
	\quad
	v^2 \in H^{10}_0, 
	\;
	v^1\in \1{H^{10}_0}.
\eeq
Here we use the subscript nation 
$q_{z_j}$ for the partial derivative,
and $\d_v$ denotes the directional derivative 
along the constant vector field
identified with the vector $v\in \{0\}\times \C^{n-1}$ in some local holomorphic coordinates.

The tensor $\t^2_p$ has the unique
$\C$-bilinear symmetric extension
$$
	\t^2_p \colon \C H_p \times \C H_p \to \C Q_p,
$$
for which we still write $\t^2_p$ by a slight abuse,
where we extend by symmetry to $\1{H^{10}_p}\times H^{10}_p$
and by zero to ${H^{10}_p}\times H^{10}_p$
and $\1{H^{10}_p}\times \1{H^{10}_p}$.
For the above quadric example, \eqref{levi-calc}
still holds for the symmetric extension.
%
The choice of the {\em symmetric $\C$-bilinear tensor}
rather than hermitian, as common for the Levi form,
will help us to keep the notation lighter for the 
subsequent
Levi form derivatives,
as that way we shall never need to remember, 
which arguments are $\C$-linear and which are $\C$-antilinear.

Recall that $M$ is pseudoconvex at $p$
if and only if there exists a nonzero covector $\theta_0\in \Om^0_p$
with
$$
\la
	\theta_p,
	\t^2_p(v, \1{v}) 
\ra
	\ge 0,
	\quad v \in H^{10}_p.
$$
We shall always assume this choice of $\theta$,
whenever $M$ is pseudoconvex.

We say that a {\em point $p\in M$ is of Levi rank $q$},
if the Levi form $\t^2_p$ at $p$ has rank $q$.
A subbundle $V\subset H^{10}$
is called {\em Levi-nondegenerate},
if the Levi form is nondegenerate on $V\times \1V$. 
For every such subbundle $V$,
we write 
$$
	V^\perp \subset H^{10},
	\quad
	V^\perp = \cup_x V_x^\perp,
	\quad
	V_x^\perp = \{ v\in H^{10}_x : \t^2_x(v, \1v^1) =0 \text{ for all } v^1\in V\},
$$
for the orthogonal complement with respect to the Levi form,
which is necessarily a subbundle.

Finally, we write $K^{10}_p\subset H^{10}_p$
and $K^{01}_p = \1{K^{10}_p}\subset H^{01}_p$
for the Levi kernel components at $p$, 
$\C K_p= K^{10}_p \otimes \1{K^{10}_p}$
for the complexification and
$K_p = \C K_p \cap T_p$
for the corresponding real part.

\bigskip
The following is an overview of some of the main results
(further results below don't assume pseudoconvexity):

\bt\Label{main}
Let $M\subset\cn$ be a pseudoconvex real hypersurface.
Then for every $q\in \{0, \ldots, n-1\}$, there exist 
an invariant submodule sheaf $\6S^{10}(q)$ of $(1,0)$ vector fields,
an invariant ideal sheaf $\6I(q)$ of complex functions,
and for every $p\in M$ of Levi rank $q$,
an invariant quartic tensor
$$
	\t^4_p \colon  \C T_p \times \C T_p 
	\times K^{10}_p\times \1{K^{10}_p} \to \C Q_p,
$$
and
a real submanifold $S\subset M$ through $p$,
such that the following hold:

\begin{enumerate}

\item 
$S$ contains the set of all points $x\in M$ of Levi rank $q$
in a neighborhood of $p$.

\item
The tangent space of $S$ at $p$ equals the real part of the kernel of $\t^4_p$:
$$
	T_p S = \Re\ker \t^4_p = \{ v\in T_p :  \t^4_p(v, v^3, v^2, v^1) = 0 
	\text{ for all } v^3, v^2, v^1 \}.
$$

\item
In suitable holomorphic coordinates vanishing at $p$,
$M$ admits the form
$$
	2\Re w = \sum_{j=1}^q |z_{2j}|^2 + \phi^4(z_4,\z_4) + o_w(1),
	\quad
	(w, z_2, z_4)\in \C\times \C^{q}\times \C^{n-q-1},
$$
where $o_w(1)$ indicates terms of weight greater than $1$,
with weights $1$, $1/2$ and $1/4$ assigned to
the components of $w$, $z_2$ and $z_4$ respectively,
and 
where $\phi^4$ is a plurisubharmonic 
homogenous polynomial of degree $4$
representing a restriction of the quartic tensor $\t^4_p$ in the sense that
$$
	\t^4_p(v^4, v^3, v^2, v^1) 
	= \d_{v^4} \d_{v^3} \d_{v^2} \d_{v^1} \phi^4
$$
holds for $v^4, v^3\in \C K_0$
and
$v^2, \1v^1\in K^{10}_0$.

\item
$S$ is given by 
$$
	S = \{f^1 = \ldots = f^m =0\},
	\quad
	df^1\wedge \ldots \wedge df^m \ne 0,
	\quad
	f^j \in \Re \6I(q).
$$
In fact, any $f\in \Re \6I(q)$ vanishes
on the set of point of Levi rank $q$. 

\item
The ideal sheaf $\6I(q)$
is generated by all functions $g$, $f$ of the form
$$
	g=\la \th, [L^2, L^1] \ra,
	\quad
	 f=L^3 \la \th, [L^2, L^1] \ra,
$$
where $\th\in\Om^0$ is a contact form, 
$L^3\in \C T$ arbitrary complex vector field,
and $L^2, \1{L^1} \in \6S^{10}(q)$
arbitrary sections of the submodule sheaf.

\item
The submodule sheaf $\6S^{10}(q)$
contains all germs of $(1,0)$ vector fields $L$
satisfying $L\in V_L^\perp$,
with $V_L\subset H^{10}$
being some Levi-nondegenerate subbundle
of rank $q$ in a neighborhood of $p$
(that may depend on $L$).
In particular, $\6S^{10}(q)$ generates
the Levi kernel at each point of Levi rank $q$.

\item
The tensor $\t^4_p$ has the positivity property
$$
	\t^4_p(v^2, v^2, v^1, \1{v^1}) \ge 0,
	\quad
	v^2 \in T_p,
	\quad
	v^1 \in K^{10}_p. 
$$

\end{enumerate}

In addition, when $M$ is of finite type at most $4$ at $p$,
the following also holds:

\begin{enumerate}

\item[(i)]
The intersection $T_p S\cap K_p$ with the Levi kernel $K_p$ is totally real.
\item[(ii)]
For every $v\in K^{10}_p$, the tensor $\t^4_p$ does not identically 
vanish on
$$
	(\C v + \C \1v) \times (\C v + \C \1v) \times \C v\times \C \1v.
$$
In particular, the regular type at $p$ equals to the 
D'Angelo type and is either $2$ or $4$.
The type is $4$ whenever $K^{10}_p\ne 0$.
\item[(iii)]
The (Catlin's) multitype at $p$ equals
$$
	(1,
		\underbrace{2,\ldots, 2}_q, 
		\underbrace{4,\ldots, 4}_{n-q-1}
	),
$$
where the number of $2$'s equals the Levi rank $q$ at $p$.
In particular, the multitype is determined by the Levi rank.

\end{enumerate}
\et

Note that we prefer the reversed order of the vector fields
$L^3, L^2, L^1$,
e.g.\
$L^3 \la \th, [L^2, L^1] \ra$
that better reflects the logical
order of the operations:
first form the Lie bracket  
$ [L^2, L^1]$
inside,
then differentiate by $L^3$
from outside (after pairing with $\th$).

For  the proofs and more detailed and general statements,
see the respective sections below.
The submodule sheaves $\6S^{10}(q)$ from Part (6)
are defined in \S\ref{submodules},
and the ideal sheaves $\6I(q)$ from Part (5)
in \S\ref{ideal-par}.
In particular, local sections in $\6I(q)$
vanish at points of Levi rank $q$
by Corollary~\ref{iq-van}.
The quartic tensor $\t^4$ is constructed 
\S\ref{quartic}.
In view of Proposition~\ref{main0},
the intersection of real kernels of differentials
$df$ for $f\in \Re \6I(q)$ coincides
with $\Re \ker \t^4_p$.
Hence we can choose functions $f^j$
satisfying (4) and (2).
The normal form in (3) 
follows from Proposition~\ref{4-normal}.

When $M$ is of finite type $4$,
Proposition~\ref{type-quartic}
implies that $\t^4$
has no holomorphic kernel,
and therefore its (real) kernel as in (2) is totally real,
as stated in (i).
Statement (ii) is also part of Proposition~\ref{type-quartic}.
Finally, statement (iii) about the multitype 
is contained in \S\ref{m-type}.

\bigskip

The following simple example
illustrates one of the differences between
functions in the ideal sheaf $\6I(q)$
and  the boundary systems
(as defined in \cite[\S2]{C84ann},
see also \S\ref{bs} below).

\be\Label{transv-f}
Consider the hypersurface
$M\subset\C^2_{w,z}$ given by
$$
	2\Re w = \phi(z,\z, \Im w), \quad 
	\phi(z,\z,u):= |z|^4 + u^2|z|^2,
$$
which is pseudoconvex and of finite type $4$.
Then a boundary system $\{ L_2; r_2\}$
%
defines the  $2$-dimensional
submanifold $S:=\{ r_2 =0\} \subset M$,
which contains all points of Levi rank $0$.
However, since $r_2$ is of the form
\beq\Label{r2}
	r_2 = \Re L^3 \la \th, [L^2, L^1]\ra
\eeq
(cf.\ the notation of Theorem~\ref{main}, part (5)),
its differential at $0$ is given by
$$
	d r_2 (v) = \Re \t^4_0 (v, L^3_0, L^2_0, L^1_0),
$$
which vanishes on the transversal space $\{dz=0\}$.
Consequently, any $S$ defined by a boundary system function $r_2$
must be tangent to the real line $\{dz=0\}$.

On the other hand, 
in the ideal sheaf $\6I(0)$
we can choose a function given by \eqref{r2}
with {\em transversal $L^3$},
i.e.\ $L^3_p\notin \C H_p$.
That allows to reduce the submanifold 
$S$ in Theorem~\ref{main} down to only the origin 
$z=w=0$,
which, in fact, is the set of points of Levi rank $0$.

In particular, the ideal sheaf $\6I(0)$ captures
as its zero set with linearly independent differentials
the $0$-dimensional singular stratum of all Levi-degenerate points,
which cannot be achieved with boundary systems.

\ee

\section{Invariant cubic tensors}

We begin by investigating the $3$rd order invariants
without the pseudoconvexity assumption.

\subsection{Double Lie brackets}
In presence of a nontrivial Levi kernel $K^{10}_p$
at a point $p\in M$,
it is natural to look for cubic tensors
arising from double Lie brackets
with one of the vector fields 
having its value inside the Levi kernel
at the reference point: 
\beq\Label{triple}
	\la \th, [L^3, [L^2, L^1]] \ra, 
	\quad
	\th \in \Om^0,
	\quad L^3, L^2\in H^{10}, 
	\quad L^1\in \1{H^{10}}, 
	\quad
	L^1_p \in \1{K^{10}_p}.	
\eeq

However, the following simple example shows that 
\eqref{triple} does not define a tensor in general:

\be\Label{one-ker}
Let $M\subset \C^3_{z_1, z_2, w}$ be
the degenerate quadric
\beq\Label{deg-qua}
	\r = -(w+\w) + z_1\z_1 =0,
\eeq
and consider the $(1,0)$ vector fields
$$
	L^3:= \d_{z_2},
	\quad
	L^2:= \d_{z_2} + cz_2 L^1,
	\quad
	L^1:= \d_{z_1} + \z_1\d_{w}.
$$
The main idea here is 
to ``twist'' the
vector field $L^2$ with $L^2_0\in K^{10}_0$
by adding a multiple of the other vector field $L^1$,
along which the Levi form is nonzero.

Then 
$$
	[L^3, [L^2, \1{L^1}]] = c [L^1, \1{L^1}] = c (\d_{\w} - \d_{w}),
$$
and hence for any fixed contact form $\th\in \Om^0$,
the value 
$$
	\la \th, [L^3, [L^2, \1{L^1}]] \ra_0
$$
depends on $c$,
even though all values $L^j_0$ are independent of $c$.
Note that both $L^2_0$ and $L^3_0$ (but not $\1{L^1_0}$)
are inside the Levi kernel $K^{10}_0$.
Hence the double Lie bracket does not define any tensor
$K^{10}_p\times K^{10}_p \times \1{H^{10}_p}\to \C Q_p$
with $p=0$.

Similarly, taking
$$
	L^3:= \d_{z_2},
	\quad
	L^2:= \d_{z_1} + \z_1\d_{w},
	\quad
	L^1:= \d_{z_2} + c\z_2 L^1,
$$
we conclude that
$$
	\la \th, [L^3, [L^2, \1{L^1}]] \ra_0
$$
depends on $c$,
thus also not defining any tensor
$K^{10}_p\times H^{10}_p \times \1{K^{10}_p}\to \C Q_p$.

The same example also shows that
the Levi form derivative $L^3 \la \th, [L^2, \1{L^1}]\ra$
considered below does not
behave tensorially on the same spaces.
\ee

In contrast, we do get an invariant tensor when {\em both vector fields
inside the inner Lie bracket have their values in the Levi kernel
at the reference point}.
We write $\t^{31}$ for the corresponding tensor, 
emphasizing the fact that it will become a restriction 
of the full tensor $\t^3$ below.

\bl\Label{bracket-tensor}
	The double Lie bracket $[L^3, [L^2, L^1]]$ defines
	an invariant tensor
	\beq
		\t^{31}_p\colon \C K_p \times K^{10}_p \times \1{K^{10}_p} \to \C Q_p,
	\eeq
	i.e.\ there exists an unique $\t^{31}_p$ as above satisfying
	$$
		\t^{31}_p(L^3_p, L^2_p, L^1_p) 
		= i  [L^3, [L^2, L^1]]_p \mod \C H_p
		,
		\quad L^3, L^2, \1{L^1} \in H^{10},
		\quad L^3_p, L^2_p, \1{L^1_p}\in K^{10}_p.
	$$
	
	Furthermore, $\t^{31}_p$ is symmetric on $K^{10}_p\times K^{10}_p\times \1{K^{10}_p}$ 
	in its $K^{10}$-arguments, 
	and on $\1{K^{10}_p}\times K^{10}_p\times \1{K^{10}_p}$ in its $\1{K^{10}}$-arguments,  
	and satisfies the reality condition
\beq\Label{30-sym}
	\1{\t^{31}_p(v^3, v^2, v^1)} = \t^{31}_p(\1v^3, \1v^1, \1v^2).
\eeq
\el
\bpf
It suffices to show that
\beq\Label{van}
	[L^3, [L^2, L^1]]_p \in \C H_p
\eeq
holds
whenever any of the values $L^j_p$ is $0$.
Since any such $L^j$ can be written as linear combination
$\sum a_k L_k$ with $a_k(p)=0$ and $L_k$ being in the same bundle
(either $H^{10}$ or $\1{H^{10}}$),
it suffices to assume $L^j = a \2L^j$ with $a(p)=0$.
Then, any term giving a nonzero value at $p$ in
 \eqref{van},
 must involve differentiation of
 the function $a$,
 either by one of the vector fields $L^j$,
 with the bracket of the other two as factor:
 $$
 	(L^{j_3}a) [L^{j_2}, L^{j_1}],
 $$
 or by two of the vector fields $L^j$,
 with the third one as factor:
 $$
 	(L^{j_3} L^{j_2}  a) L^{j_1},
 $$
In the second case \eqref{van} is clear.
In the first case, we obtain
a bracket of two $L^j$,
one of which has value at $p$ contained in
the kernel $K^{10}_p\oplus \1{K^{10}_p}$,
implying \eqref{van}.

The reality condition is straightforward and the symmetries follow
from the Jacobi identity.
\epf

\br
A closely related construction 
is the one of the cubic form
$c=c(L^3, L^2, L^1)$
 by Webster~\cite{W95},
defined for triples of vector fields in a neighborhood
of a reference point $p\in M$,
whose value $c_p$ at $p$
depends on the $1$-jets of the vector fields $L^3, L^2, L^1$,
but in general, is not uniquely determined by their values at $p$,
as demonstrated by Example~\ref{one-ker},
unless all three vector fields are valued in kernel at $p$
as in Lemma~\ref{bracket-tensor}.
\er

\subsection{The Levi form derivative}
As alternative to the double Lie bracket tensor, one can 
differentiate the Levi form after pairing with a contact form,
which is similar to the approach employed by Catlin 
in his boundary system construction:
\beq\Label{der}
	L^3 \la \th, [L^2, L^1] \ra.
\eeq
Again Example~\ref{one-ker} shows that \eqref{der}
does not define a tensor if either
of the vector fields $L^2$, $L^3$ inside the Lie bracket
is not in the Levi kernel at $p$.
On the other hand, if {\em both vector fields
$L^1$, $L^2$ inside the Lie bracket
have their value at $p$ contained in the Levi kernel}, we do obtain a tensor
{\em even when the outside vector field $L^3$ 
is not necessarily contained in $\C H$}:

\bl\Label{levi-der}
There exists unique cubic tensor
\[
	\t^3_p \colon \C T_p \times  K^{10}_p \times \1{K^{10}_p} \to \C Q_p,
\]
satisfying
	$$
		\la \th_p, \t^3_p(L^3_p, L^2_p, L^1_p) \ra
		= i (L^3\la \th, [L^2, L^1] \ra)_p
		,
		\quad \th\in \Om^0,
		\quad L^3 \in \C T,
		\quad L^2, \1{L^1} \in H^{10},
		\quad L^2_p, \1{L^1_p}\in K^{10}_p.
	$$

Furthermore, $\t^3_p$ satisfies the reality condition
\beq\Label{3-sym}
\1{\t^3_p(v^3, v^2, v^1)} = \t^3_p(\1{v^3}, \1{v^1}, \1{v^2}),
\eeq
note the switch of the last two arguments.
\el

\bpf
The proof is similar to that of Lemma~\ref{bracket-tensor},
and the symmetry follows directly from the definition.
\epf

\br
Note that the tensor $\t^3$ in Lemma~\ref{levi-der}
is defined when the first argument is {\em any complex vector field} 
on $M$, not necessarily in the subbundle $\C H$,
in contrast to the tensor $\t^{31}$ in Lemma~\ref{bracket-tensor}.
This shows that taking derivatives of the Levi form
provides more information than taking iterated Lie brackets.
\er

\subsection{A normal form of order $3$ and 
complete set of cubic invariants of a real hypersurface}
\Label{norm-f}
To compare tensors $\t^{31}$ and $\t^3$,
it is convenient to use a partial normal form
for the cubic terms.
In the following we write $\phi^{j_1\ldots j_m}(x_{j_1},\ldots x_{j_m})$
for a polynomial of the multi-degree $(j_1,\ldots, j_m)$
in its corresponding variables.
We also write 
$z_k=(z_{k1},\ldots, z_{km})\in \C^m$
for the coordinate vectors and their components.

\bp\Label{3-normal}
For every real hypersurface $M$ in $\C^{n}$ and point $p\in M$ of Levi rank $q$, 
there exist local holomorphic coordinates 
$$
	(w,z)=(w, z_2, z_3)\in \C\times \C^q\times \C^{n-q-1},
$$ 
vanishing at $p$,
where $M$ takes the form
$$
	w + \w = \phi(z,\z, i(w-\w)), 
	\quad
	\phi(z,\z,u) = \phi^2(z,\z,u) + \phi^3(z,\z,u) + O(4), 
$$
where
$$
	\phi^2(z,\z,u) = \phi^{11}(z_2, \z_2) = \sum \pm |z_{2j}|^2,
$$
and
$$
	\phi^3(z,\z,u) = 2\Re \phi^{21}(z,\z_3) + \phi^{111}(z_3, \z_3, u),
$$
with
$$
	\phi^{21}=
	\sum\phi^{21}_{j_1j_2j_3} z_{2j_1} z_{2j_2} \z_{3j_3}
	+
	\sum\phi^{21}_{j_1j_2j_3} z_{2j_1} z_{3j_2} \z_{3j_3}
	+
	\sum\phi^{21}_{j_1j_2j_3} z_{3j_1} z_{3j_2} \z_{3j_3}
	,
$$
$$
	\phi^{111} = \sum \eps_j |z_{3j}|^2 u,
	\quad \eps_j\in \{-1, 0, 1\},
$$
and $O(4)$ stands for all terms of total order at least $4$.
\ep

\bpf
It is well-known that the quadratic term $\phi^2$ 
can be transformed into $\phi^{11}(z^2, \z^2)$ 
representing the (nondegenerate part of) the Levi form.
Furthermore, as customary, we may assume that the cubic term $\phi^3$ has no harmonic terms.

Next, by suitable polynomial transformations 
$$
	(z,w)\mapsto (z + \sum_{j=1}^r z^2_j h_j(z,w) , w),
$$
we can eliminate all cubic monomials of the form $\z^2_j h(z,u)$
and their conjugates, where $h(z,w)$ is any holomorphic quadratic monomial.
The proof is completed by inspecting the remaining cubic monomials
and diagonlizing the quadratic form in the component $\phi^{111}$.
\epf

Next we use the convenient $(1,0)$ vector fields with obvious notation:

\bl\Label{vf-norm}
For a real hypersurface $M\subset \C^{n}$ given by
\beq\Label{graph}
	w + \w = \phi(z,\z, i(w-\w)), \quad (w, z)\in\C\times \C^{n-1},
\eeq
the subbundle $H^{10}$ of $(1,0)$ vector fields is spanned by
$$
	L_j := \d_{z_j} + \frac{\phi_{z_j}}{1-i\phi_u} \d_w, 
	\quad 
	j=1,\ldots, n-1.
$$
More generally, $H^{10}$ is spanned by all vector fields of the form
\beq\Label{lv}
	L_v:= \d_v + \frac{\phi_{v}}{1-i\phi_u} \d_w, 
	\quad v\in \{0\} \times \C^{n-1},
\eeq
where the subscript $v$ denotes the differentiation in the direction of $v$.
\el

Calculating with 
the special vector fields from Lemma~\ref{vf-norm}, we obtain:
\bc\Label{3-calc}
Let $M$ be in the normal form given by Proposition~\ref{3-normal}.
Then tensors $\t^{31}_p$ and $\t^3_p$ defined in Lemmas~\ref{bracket-tensor}
and ~\ref{levi-der} respectively satisfy
\beq\Label{3-dif}
	\la \th_0,  \t^{31}_0(v^3, v^2, v^1) \ra
	=
	\la \th_0,  \t^3_0(v^3, v^2, v^1) \ra
	= \d_{v^3} \d_{v^2} \d_{v^1} \phi^3,
\eeq
where
$$
	v^3, v^2, \1{v^1} \in K^{10}_0 \cong \{0\} \times \C^{n-q-1},
	\quad
	\th=i \d\rho, 
	\quad
	\rho =-2\Re w + \phi.
$$
Furthermore, the second identity in \eqref{3-dif} still holds for $v^3\in \C H_0$.
\ec
In particular,
$\t^{31}$ is a restriction of $\t^3$ 
to $\C T_p \times  K^{10}_p \times \1{K^{10}_p}$,
explaining the notation.

\br
The term $\phi^{21}$ in Proposition~\ref{3-normal}
represents, up to a nonzero constant multiple,
the cubic invariant tensor 
introduced by Ebenfelt~\cite{E-jdg},
defined
by means of the Lie derivative $\6T$:
\beq\Label{lie}
	\t^{21}_p\colon 
	H^{10}_p\times H^{10}_p \times \1{K^{10}_p}
	\to \C Q_p,
	\quad
	\la \th_p, \t^{21}(L^3_p, L^2_p, L^1_p) \ra =
	\la \6T_{L^3} \6T_{L^2} \th, L^1 \ra_p,
\eeq
where
$$
	\6T_L := d \circ \imath_L + \imath_L \circ d,
$$
and $\imath$ is the contraction.
Here the phenomenon illustrated by
 Example~\ref{one-ker}
of the lack of tensoriality in the last argument
does not occur as 
the right-hand side obviously depends only
on the value $L^1_p$.
In fact, it follows from the
tranformation law of the Lie derivative,
$$
	\6T_L(f\th) = f\6T_L\th + (Lf)\th ,
	\quad
	\6T_{fL}\th = f \6T_L\th + \th(L) df,
$$
that the same right-hand side in \eqref{lie}
defines a tensor even on the larger spaces
\beq\Label{21}
	\t^{21}_p\colon 
	\C H_p\times \C K_p \times \C K_p\to \C Q_p,
\eeq
that we denote by the same letter in a slight abuse of notation.

On the other hand, the same expression 
{\em does not 
define any tensor when the first argument 
varies arbitrarily
in $\C T$}, even for a Levi-flat hypersurface $M\subset\C^2$.
Indeed, taking 
$$
	M=\{ \rho = 0\},
	\quad
	\rho= -w - \w,
$$ 
	with
$$
	\th = -i (1 + z+\z) d w,
	\quad
	L^3 = i(\d_w-\d_\w),
	\quad
	L^2 = \d_{z},
	\quad
	L^1 = \d_\z,
$$
we compute
$$
	\6T_{fL^3} \6T_{L^2} \th 
	=\6T_{fL^3} (\imath_{L^2} d\th)
	= f \omega +  d\th(L^2, L^3) df
	= f\omega + df,
$$
where $f$ is any smooth complex function and
 $\omega$ is some $1$-form.
 Then choosing $f$ with $f( p) =0$
 and $df_p=cd\z$, we conclude that 
 $\la \6T_{fL^3} \6T_{L^2} \th, L^1\ra_p$
 depends on $c$ and hence is not tensorial.
\er

Finally, as consequence from Proposition~\ref{3-normal},
we obtain:
\bc
The tensors
 $\t^3$ and $\t^{21}$,
 given by respectively Lemma~\ref{levi-der} and \eqref{21},
coincide up to a constant on their common set of definition,
and
constitute together the {\em full set of cubic invariants} of $M$ at $p$.
\ec

\subsection{Symmetric extensions}
As consequence Corollary~\ref{3-calc},
$\t^3$ is symmetric in
$K^{10}$- or in $\1{K^{10}}$-vectors
whenever two of them occur in any two arguments.
This property leads to a natural symmetric extension:

\bl\Label{3-symm}
The restriction 
$$
	\t^{30}_p \colon \C K_p \times K^{10}_p\times \1{K^{10}} \to \C Q_p
$$ 
of the cubic tensor $\t^3_p$
 admits an unique symmetric extension
$$
	\2\t^{30}_p \colon
	\C K_p \times \C K_p \times \C K_p \to \C Q_p,
$$
satisfying
$$
	\la \th_0,  \2\t^{30}_0(v^3, v^2, v^1) \ra
	= \d_{v^3} \d_{v^2} \d_{v^1} \phi^3,
$$
whenever $M$ is in a normal form $\rho= -2\Re w +\phi=0$
as in Proposition~\ref{3-normal}
and $\th = i\d\rho$.
\el

\br
Note that since $\phi^3$ has no harmonic terms in a normal form,
the extension tensor $\2\t^{30}$ 
always vanishes whenever its arguments
are either all in $K^{10}$ or all in $\1{K^{10}}$.
\er

\be
In contrast to $\t^{30}$, the full cubic tensor $\t^3$
does not in general have any invariant extension to
$\C T\times \C K \times \C K$.
Indeed, consider the cubic $M\subset \C^2$ given by
$$
	\rho := -2\Re w + \phi ^3=0,
	\quad
	\phi^3 = 2\Re (z^2\z).
$$
Then $\d_{\w} \d_z \d_z\phi^3=0$.
Now consider a change of coordinates with linear part
$(w,z)\mapsto (w, z+ iw)$ transforming $\d_{\w}$
into $\d_{\w}-i\d_{\z}$.
Then, after removing harmonic terms, the new cubic term takes form
$$
	\phi^3 = 2\Re (z^2\z) - 4\Im w z\z.
$$
But then $(\d_{\w} -i\d_{\z}) \d_z \d_z\phi^3\ne 0$, 
i.e.\ the $3$rd derivatives of $\phi^3$ do not transform as tensor
when passing to another normal form.
\ee

%
%


%

\subsection{Cubic tensors vanishing for pseudoconvex hypersurfaces}
If $M$ is pseudoconvex, the Levi form $\la \th, [L, \1{L}]\ra$
does not change sign, and therefore the cubic tensor $\t^3$ must vanish identically. We obtain:

\bl\Label{psc-vanish}
Let $M$ be a pseudoconvex hypersurface and $p\in M$.
Then the cubic tensor $\t^3_p$ (and therefore its restriction $\t^{31}_p$) vanishes identically.
Equivalently, the cubic normal form in Proposition~\ref{3-normal} satisfies
\beq\Label{psc-cubic}
	\phi^{21}(z, \z_3) = \sum_{jkl} c_{jkl} z_{2j} z_{2k} \z_{3l},
	\quad
	\phi^{111}(z_3, \z_3, u) =0.
\eeq
\el

The remaining cubic terms in \eqref{psc-cubic} 
can be absorbed into higher weight terms as follows.
We write $o_w(m)$ for terms of weights heigher than $m$.

\bc\Label{pse-nf-cor}
A pseudoconvex hypersurface $M$ in suitable holomorphic coordinates
is given by
\beq\Label{psc-cubic-red}
w + \w = \phi(z,\z, i(w-\w)), 
\quad
\phi =
 \sum_j |z_{2j}|^2 + o_w(1),
\eeq
where $o_w$ is calculated for $(w, z_2, z_3)$, and their conjugates, 
having weights 
$1$, $\frac12$, $\frac13$ respectively.
In particular, $M$ cannot be of (D'Angelo) type $3$.
\ec

The last statement follows directly from \eqref{psc-cubic-red},
since the contact orders  with lines in the directions of
$z^3$
are at least $4$.

\subsection{Freeman's modules and uniformly Levi-degenerate hypersurfaces}
Freeman~\cite{Fr77} introduced for any smooth real hypersurface
$M\subset\C^n$, a decreasing sequence
of invariantly defined submodules
of the module of all smooth $(1,0)$ vector fields on $M$.
In particular, Freeman's second submodule
$N'_2$ is defined by the Lie bracket relation
$$
	N'_2 := \{ L \in H^{10} : [L, \1{L^1}]\in \C H \text{ for all } 
	L^1\in H^{10} \},
$$
or equvalently, by the inclusion relation in the Levi kernels $K^{10}_p$:
$$
	N'_2  :=  \{ L \in H^{10} : L_p\in K^{10}_p \text{ for all } p \}.
$$

For a fixed $p_0\in M$,
it is easy to see that
vector fields in $N'_2$ span
the Levi kernel $K^{10}_{p_0}$
if and only if $\dim K^{10}_p$ 
is constant for $p\in M$ near $p_0$,
i.e.\ when $M$ is {\em uniformly Levi-degenerate} 
in a neighborhood of $p_0$.
Indeed, since
$\dim K^{10}_p$ is 
upper semi-continuous,
whereas the dimension of the span
at $p$ of all values of $N'_2$
is lower semi-continuous,
the only way both dimensions can match
at $p_0$
is when they are both constant
in its neighborhood,
i.e.\ when $M$ is uniformly Levi-degenerate.
In the latter case, 
the tensor $\t^{21}$
in \eqref{lie} 
can be computed by means of
double Lie brackets of vector fields in $N'_2$,
as shown in \cite[Appendix]{KZ06}.

On the other hand, if $M$ is of finite type,
one necessarily has $K_p=0$ on a dense set of $p\in M$
(otherwise nontrivial integral surfaces of $K^{10}$
would be complex-analytc subsets of $M$).
Hence in this case, the module 
$N'_2$ is always trivial,
whereas the tensor $\t^3$ may not be so.
And even when the module $N'_2$ is not trivial
but $\dim K_p$ is not constant 
in any neighborhood of a point $p_0$,
it is easy to see 
that the set of values $L_p$ for $L\in N'_2$
can never span the full Levi kernel $K_p$

\br
In the {\em  uniformly Lev-degenerate case},
i.e.\
when the Levi kernel dimension $\dim K_p$ is constant,
alternatively to the Lie derivative approach in \eqref{lie},
both double Lie brackets (as in Lemma~\ref{bracket-tensor}) 
and Levi form derivative approaches
(as in Lemma~\ref{levi-der})
can be used to define $\t^{21}$
by 
{\em imposing additional restrictions}
on the vector fields
to be contained in the 
Levi kernel subbundle everywhere,
rather than only at the reference point
as in \eqref{lie}.

On the other hand,
without the uniformity assumption on $\dim K_p$,
only the Lie derivative approach leads
to an invariant definition of $\t^{21}$,
whereas only the Levi form derivative approach is suitable 
to define the full cubic tensor $\t^3$ as in Lemma~\ref{levi-der}.
It is quite remarkable that no single approach
seems to work to define
the complete system of cubic invariants,
consisting of the pair $(\t^{21}, \t^3)$.
\er

\section{Invariant quartic tensors}
If the cubic tensor $\t^3$ vanishes,
it is natural to look for higher order invariants by taking iterated Lie brackets
or higher order derivatives of Levi form.
However, 
in contrast to the statements of Lemmas~\ref{bracket-tensor} 
and \ref{levi-der},
we don't obtain any tensor in this way
even when all vector field values are in the Levi 
kernel,
as demonstrated by our next counter-example:

\be\Label{ex-4}
Let $M\subset \C^3_{z_1, z_2, w}$ be again
the degenerate quadric 
from Example~\ref{one-ker},
and set 
$$
	L:=\d_{z_2} +  c z_2(\d_{z_1} + \z_1\d_w).
$$
Then 
$$
	[L, \1L] = |c z_2|^2 (\d_{\w} - \d_w),
$$
and
both $\la \th, [L, [\1L, [L, \1L]]] \ra_0$
and $(L\bar L \la \th,  [L, \bar L] \ra)_0 $
depends on $c$ (and hence on the $1$-jet of $L$), even though
the value $L_0$ is contained in the Levi kernel $K_0^{10}$,
and the cubic tensor $\t^3_0$ identically vanishes.
\ee

%
%

\subsection{Vector fields that are in the Levi kernel up to order $1$}
In view of Example~\ref{ex-4},
in order to obtain a tensor, we need to restrict
the choice of the vector fields.
This motivates the following definition:

\bd\Label{ker-1}
Let $L$ be a $(1,0)$ vector field.
We say that 
{\em $L$ is in the Levi kernel up to
order $1$ at $p$} if,
for any vector fields $L^1\in H^{10}$, $L^2\in \C T$,
and any contact form $\th$, the following holds:
\beq\Label{1-ker}
	\la \th, [L^1, \1L]\ra_p
	= (L^2 \la \th, [L^1, \1L] \ra)_p 
	= 0.
\eeq

More generally, 
we have the following
``microlocal'' version of this definition as follows.
For a fixed tangent vector $v\in \C T_p$,
we say that
{\em $L$ is in the Levi kernel up to 
$v$-order $1$ at $p$} 
if \eqref{1-ker}
holds whenever $L^2_p=v$,
i.e.\ we differentiate the Levi form only in the fixed given direction.
(The latter property obviously depends only on the value $v$
rather than its vector field extension $L^2$.)
If the above property holds for all $v$
in a vector subspace $V\subset \C T_p$,
we also say that
{\em $L$ is in the Levi kernel up to 
$V$-order $1$ at $p$}.

If $L$ is a $(0,1)$ vector field,
we say that $L$ is in the Levi kernel up to order $1$ at $p$
whenever $\1L$ is.
Similarly, we extend all the other terminology in this definition
to $(0,1)$ vector fields.
\ed


It is straightforward to see that:
\bl\Label{1-jet-dep}
For any $(0,1)$ vector field $\1L$
with $\1L_p\in \1{K^{10}_p}$, 
the expression
$$ 
	(L^2 \la \th, [L^1, \1L] \ra)_p 
$$
only depends on the values $L^2_p, L^1_p$ and $\th_p$,
as well as on the $1$-jet of $L$ at $p$.
In particular,  $L$ being in the Levi kernel up to order $1$,
is a linear condition on the $1$-jet of $L$ at $p$.
\el

\be
In the setting of Example~\ref{ex-4},
choosing $L^1:= \d_{z_1} + \z_1\d_w$, we compute
$$
	(L \la \th, [L^1, \1L] \ra )_0 \ne 0,
$$
which shows that here $L$ is not in the Levi kernel of order $1$,
even though its value at $0$ is contained in the Levi kernel.
\ee

\br\Label{bracket-restrict}
Using any normal form as in Proposition~\ref{3-normal} 
and calculating with vector fields \eqref{lv},
we can obtain a condition equivalent to
 \eqref{1-ker}
with $L^2$ in $\C H$ (rather than in $\C T$), which can be 
stated in terms of the double Lie brackets instead of
the Levi form derivatives:
\beq\Label{1-ker'}
	\la \th_p, [L^1, \1L]_p\ra
	= \la \th_p, [L^2, [L^1, \1L]]_p \ra 
	= \la \th_p, [\1{L^2} , [L^1, \1L]]_p\ra)
	= 0.
\eeq
\er

A priori, it is not at all clear
that vector fields as in Definition~\ref{ker-1} exist.
The following lemma
provides an easy way of constructing them.

\bl\Label{1-ker-def}
Let $M$ have the Levi rank $q$ at $p\in M$,
with Levi kernel $K^{10}_p$.
Assume that $L\in H^{10}$ is in the Levi kernel
up to order $1$ at $p$,
as per Definition~\ref{ker-1}.
Then $L_p\in K^{10}_p$ and
\beq\Label{3-vanish}
	\t^3_p(L^2_p, L^1_p, \1L_p) 
	= 0, 
	\quad L^2\in \C T, 
	\quad
	L^1\in H^{10}, 
	\quad
	L^1_p\in K^{10}_p.
\eeq
must hold for all $L^2$, $L^1$.
(Equivalently, $\1L_p$ is contained in the kernel of $\t^3_p$ in the last argument).

Vice versa, assume that $L_p\in K^{10}_p$ and 
\eqref{3-vanish} holds.
Let
$$(\2L^1, \ldots, \2L^q)$$
be a Levi-nondegnerate system of 
$(1,0)$ vector fields at $p$
(i.e.\ the matrix $\t^2_p(\2L^j_p, \2L^k_p)$ is nondegenerate),
such that $L$ is Levi-orthogonal to each $\2L^j$, $j=1,\ldots,q$,
in a neighborhood of $p$.
Then $L$ is in the Levi kernel up to order $1$ at $p$.
\el

\bpf
The first part follows directly from the definitions.

Vice versa, since the Levi form has rank $r$ in $p$,
and $L_p$ is Levi-orthogonal to each $\2L^j$,
it follows that $L_p$ is in the Levi kernel, i.e.\
the first expression in \eqref{1-ker} must vanish.

Next, \eqref{3-vanish}
implies that the second expression in \eqref{1-ker}
vanishes whenever $L^1_p\in K^{10}_p$.
Similarly, in view of the symmetry \eqref{3-sym},
also the third expression vanishes
under the same assumption.

Finally, to a general $L^1$ with $L^1_p\notin K^{10}_p$,
we can always add a linear combination of $\2L^j$
to achieve the inclusion of the value at $p$ in the Levi kernel.
Since $L$ is Levi-orthogonal to each $\2L^j$ identically in a neighborhood of $p$,
this does not change \eqref{1-ker},
completing the proof.
\epf

In particular, in view of Lemma~\ref{psc-vanish} we obtain:

\bc\Label{gen-pt}
Let $M$ be pseudoconvex.
Then every $v\in K^{10}_p$
extends to a $(1,0)$ vector field,
which is in the Levi kernel up to order $1$ at $p$.
\ec

\br\Label{kernel-span}
More generally, a similar result
can be obtained without pseudoconvexity
for a $v\in K^{10}_p$ 
whose conjugate $\1v$ is in the kernel of $\t^3$
(in the last argument), 
i.e.\ satisfying 
$$
	\t^3_p(L^3_p, L^2_p, \1v) = 0
$$
for all $L^3$, $L^2$.
Then there exists a $(0,1)$ vector field $\1L$ extension of  $\1v$,
which is in the Levi kernel up to order $1$ at $p$.
\er

\subsection{Invariant submodule sheaves of vector fields}
\Label{submodules}
The notion of the {\em Levi kernel inclusion up to order $1$} 
has been defined pointwise in Definition~\ref{ker-1}.
In order to have a uniform control for 
Levi kernels in nearby points,
we shall need to define
corresponding sheafs of
submodules of vector fields as follows. 

\bd\Label{q-sheaf}
Let $M\subset \C^n$ be a real hypersurface.
Denote by $\6T^{10}$ the sheaf of all 
 $(1,0)$ vector fields on $M$.
For every $q\le n-1$,
define $\6S^{10}(q)\subset \6T^{10}$ to be the submodule sheaf
consisting of all germs of vector fields on $M$ which are
contained in the Levi kernel up to order $1$
 at every point of Levi rank $\le q$.
\ed

Clearly we have the inclusions
$$
	\6S^{10}(0) \supset \cdots \supset \6S^{10}(n-1).
$$

As a direct consequence of Lemma~\ref{1-ker-def},
we obtain the following strengthening of Corollary~\ref{gen-pt}:
\bc\Label{q-sheaf-sect}
Let $M$ be a pseudoconvex hypersurface.
Then for every $q$, local sections of ${\6S}^{10}(q)$
span the Levi kernel $K^{10}_x$
at every point $x\in M$ of Levi rank $q$.
\ec

Note that to guarantee the existence of
sufficiently many sections as in 
Corollary~\ref{q-sheaf-sect},
 it is important to restrict
the property underlying Definition~\ref{q-sheaf}
only to points of Levi rank $\le q$.
Without that restriction, 
the sheaf would become trivial
e.g.\ for any
manifold $M$ that is generically 
Levi-nondegenerate
(which is the case for any $M$ of finite type).

Definition~\ref{q-sheaf}
requires to check the condition at every point of Levi rank $\le q$,
which can be difficult to deal with in practice,
when the set of such points is not ``nice''.
However, Lemma~\ref{1-ker-def} implies:

\bc\Label{span-ker}
Suppose that $M$ is pseudoconvex.
Let $V\subset H^{10}$ be
a Levi-nondegenerate subbundle
of rank $q$ in a neighborhood of $p$.
Then any section in the Levi-orthogonal complement
$V^\perp$
is contained in the sheaf $\6S^{10}(q)$.
In particular, local sections of $\6S^{10}(q)$
span the Levi kernel $K^{10}_x$
at every point $x\in M$ of Levi rank $\le q$.
\ec

Recall that in Section~\ref{results},
we call a subbundle $V\subset H^{10}$ 
 Levi-nondegenerate
whenever the Levi form restriction
to $V\times \1V$ is nondegenerate.

%

\be\Label{ex-module-0}
If $p\in M$ is a point of Levi rank $0$,
where the cubic tensor $\t^3_p$ vanishes,
any $(1,0)$ vector field is automatically
contained in the Levi kernel up to order $1$ at $p$
in view of Corollary~\ref{gen-pt}.
In particular, if $M$ is pseudoconvex,
the sheaf ${\6S}^{10}(0)$
consists of all germs of $(1,0)$ vector fields.
When $M$ is not pseudoconvex,
the condition $f\in {\6S}^{10}(0)$
is more delicate,
requiring $f$ to belong to the kernel
of the cubic tensor $\t^3_p$ at every point of Levi rank $0$.

On the opposite end, for $q=n-1$,
the sheaf ${\6S}^{10}(n-1)$ 
consists of all germs of $(1,0)$ vector fields 
that are contained in the Levi kernel at every point.
Indeed, any point is of Levi rank $\le n-1$,
hence any germ in  ${\6S}^{10}(n-1)$ 
is automatically contained in the Levi kernel at every point.
Vice versa, for any such germ $L$,
the first term in \eqref{ker-1} vanishes identically
and hence also the second vanishes by differentiation.
\ee

The sheaves  ${\6S}^{10}(q)$ for $1<q<n-1$ are more interesting: 

\be\Label{z12}
Let $M\subset\C^3_{z_1, z_2, w}$ be given by
$$
	2\Re w = |z_1|^4 + |z_2|^4.
$$
The condition $L\in {\6S}^{10}(1)$
only involves points of Levi rank $\le 1$,
i.e.\ the subset  $M^1:=\{z_1z_2=0\}$.
Hence, at points outside $M^1$, the sheaf ${\6S}^{10}(1)$
contains all germs of $(1,0)$ vector fields.

Next, for $p=(z_1, 0, w) \in M^1$ with $z_1\ne 0$,
a $(1,0)$ vector field $L$ satisfying
\beq\Label{a12}
	L = a^1 \d_{z_1} + a^2 \d_{z_2} \mod (\d_{\w} - \d_w)
\eeq
for some functions $a^1, a^2$,
is in the Levi kernel up to order $1$ at $p$
if and only if
the coefficient $a^1$ vanishes up to order $1$ at $p$.
A similar property holds for $p=(0, z_2, w)\in M^1$.
Finally, for $p=(0,0,w)\in M^1$, any $L$ is 
in the Levi kernel up to order $1$ at $p$.
However, any $L\in {\6S}^{10}(1)$
must have both $a^1, a^2$ vanish at $0$
up to order $1$ by continuity.
Summarizing, a germ $L\in {\6S}^{10}(1)$
if and only if $a^j$ vanishes up to order $1$ 
at every point of $M^1$ with $z_{3-j}=0$ for $j=1,2$.
\ee


\subsection{Construction of the quartic tensor}
\Label{quartic}
Equipped with special vector fields as in Definition~\ref{ker-1},
we can now define an invariant quartic tensor
by means of the second order derivatives of the Levi form:

\bl\Label{d2}
Let $M$ be such that the cubic tensor $\t^3_p$ vanishes for some $p\in M$.
Then there exists an unique tensor
$$
	\t^4_p \colon  \C T_p \times \C T_p 
	\times K^{10}_p\times \1{K^{10}_p} \to \C Q_p,
$$
such that for any $(1,0)$ vector fields
$\1{L^1}, L^2\in H^{10}$ that are in the Levi kernel
up to order $1$ at $p$,
any vector fields $L^3, L^4\in \C T$,
and any contact form $\th\in \Om^0$,
\beq\Label{1-levi'}
	\la \th_p, \t^4_p(L^4_p, L^3_p, L^2_p, L^1_p)\ra 
	= i(L^4 L^3 \la \th, [L^2, L^1]\ra)_p.
\eeq

More generally, \eqref{1-levi'} still holds
whenever both $L^1$ and $L^2$
are in the Levi kernel
up to $L^j_p$-order $1$ at $p$, for $j=3,4$.
\el

\bpf
Similar to the proof of Lemma~\ref{bracket-tensor},
it suffices to prove that the right-hand side of \eqref{1-levi'}
vanishes whenever
either $L^k=a\2L^k$ for some $k=1,2,3,4$, 
or $\th = a\2\th$, where $a$ is a smooth function vanishing at $p$.
In the following $\2a$ will denote either $a$ or the conjugate $\1a$
and we assume (without loss of generality)
that each of $L^3, L^4$
is contained in either $H^{10}$ or $H^{01}$.

Now the vanishing of the right-hand side in \eqref{1-levi'} 
is obvious for $k=4$.
For $k=3$, it takes
the form
$$
	(L^4\2a)_p (\2L^3 \la \th, [L^2, L^1] \ra)_p,	
$$
which must vanish in view of Definition~\ref{ker-1}.
For $k=1$, we obtain
$$
	(L^4 \bar a)_p (L^3 \la \th, [L^2, \2L^1] \ra)_p
	+ 	(L^3 \bar a)_p (L^4 \la \th, [L^2, \2L^1] \ra)_p
	+ 	(L^4 L^3 \bar a)_p (\la \th, [L^2, \2L^1] \ra)_p
	,
$$
which again vanishes in view of Definition~\ref{ker-1}.
For $k=2$, the proof follows from the case $k=1$
by exchanging $L^2$ and $L^1$ and conjugating.
Finally, for $\th = a \2\th$, we obtain
$$
	(L^4 a)_p (L^3 \la \2\th, [L^2, L^1] \ra)_p
	+ 	(L^3 a)_p (L^4 \la \2\th, [L^2, L^1] \ra)_p
	+ 	(L^4 L^3 a)_p (\la \2\th, [L^2, L^1] \ra)_p,
$$
which vanishes by the same argument.
\epf

\br
In higher generality, when the cubic tensor $\t^3_p$ may not vanish completely,
a quartic tensor $\t^4_p$ can still be constructed via \eqref{1-levi'} 
along certain kernels of $\t^3_p$.
We will not pursue this direction as our focus here is on 
the pseudoconvex case when $\t^3_p$ always vanishes identically.
\er

\subsection{Positivity of the quartic tensor}
As direct consequence of Lemma~\ref{d2} we obtain:
\bc
Let $M$ be pseudoconvex.
Then the quartic tensor $\t^4_p$ satisfies the following positivity property:
$$
	\t^4_p(v^2, v^2, v^1, \1{v^1}) \ge 0,
	\quad
	v^2 \in  T_p,
	\quad
	v^1 \in K^{10}_p. 
$$
\ec

\bpf
Since the Levi form $\t^2_p(v^1, \1{v^1})$ vanishes for all 
$v^1\in K^{10}_p$,
the function $x\mapsto i\la \th, [L^2, L^1]\ra_x$ in \eqref{1-levi'}
for fixed $L^2$ and $L^1 = \1{L^2}$
achieves its local minimum at $p$.
Hence its differential also vanishes at $p$ and
the real hessian is positive semidefinite.
Then the desired conclusion follows from \eqref{1-levi'}.
\epf

\subsection{A normal form up to weight $1/4$}
Since the cubic normal form for pseudoconvex hypersurfaces \eqref{psc-cubic-red}
is in some sense ``lacking nondegenerate terms'', 
we extend it by lowering the weight of $z_3$ from $1/3$ to $1/4$
(and renaming $z_3$ to $z_4$ to reflect the weight change):

\bp\Label{4-normal}
For every pseudoconvex real hypersurface $M$ in $\C^{n}$ 
and point $p\in M$ of Levi rank $q$, 
there exist local holomorphic coordinates 
$$
	(w,z)=(w,z_2, z_4)\in \C\times \C^q\times \C^{n-q-1},
$$ 
vanishing at $p$,
such that $M$ takes the form
\beq\Label{phi24}
	\rho = 0,
	\quad
	\rho = -2\Re w +  \phi(z,\z , i(w-\w)), 
	\quad
	\phi = \phi^2 + \phi^4 + o_w(1), 
\eeq
where
$$
	\phi^2(z,\z,u) = \sum_{j=1}^q |z_{2j}|^2,
	\quad
	\phi^4(z,\z,u) = 2\Re \phi^{31}(z_4,\z_4) + \phi^{22}(z_4,\z_4),
$$
with
$$
	\phi^{31} = 
	\sum \phi^{31}_{j_1 \ldots j_4} z_{4j_1} z_{4j_2} z_{4j_3} \z_{4j_4}
	\quad
	\phi^{22} =
	\sum \phi^{22}_{j_1 \ldots j_4} z_{4j_1} z_{4j_2} \z_{4j_3} \z_{4j_4},
$$
where the weight estimate $o_w$ is calculated for $u$, $z^2_j$, $z^4_k$, and their conjugates, 
being assigned the weights 
$1$, $\frac12$, $\frac14$ respectively.
Each polynomial $\phi^{jk}$ here is bihomogenous of bidegree $(j,k)$ in 
$(z^4, \z^4)$.

Furthermore, the following hold:
\begin{enumerate}
\item
For every 
$v\in K^{10}_0\cong  \{0\} \times \C^{n-q-1}$, 
the vector field $L_v$ given by \eqref{lv} is 
in the Levi kernel up to $v^0$-order $1$ at $0$
for any $v^0\in \C K_0$.

\item
For $v^4, v^3\in \C K_0$
and
$v^2, \1v^1\in K^{10}_0$,
we have
\beq\Label{phi-diff}
	\t^4_p(v^4, v^3, v^2, v^1) 
	= \d_{v^4} \d_{v^3} \d_{v^2} \d_{v^1} \phi^4.
\eeq
%
In particular, the restriction
\beq\Label{t40}
	\t^{40}_p \colon \C K_p\times  \C K_p \times K^{10}_p\times \1{K^{10}_p} \to \C Q_p
\eeq 
of $\t^4_p$ is symmetric in whatever arguments can be exchanged
and satisfies the reality condition
$$
	\1{\t^4_p(v^4, v^3, v^2, v^1)}
	=
	\t^4_p(\1{v^4}, \1{v^3}, \1{v^1}, \1{v^2}),
$$
note the switch of the last two arguments.
\end{enumerate}
\ep

\bpf
The existence of the desired normal form 
is a direct consequence of Lemma~\ref{psc-vanish}.
A direct calculation shows the special vector fields in \eqref{lv}
with $v\in \{0\} \times \C^{n-q-1}$ are in the Levi kernel up to tangential order $1$
as claimed.
The remaining properties are straightforward. 
\epf

Similarly to Corollary~\ref{3-calc},
one can also show the following quartic Lie bracket
representation:
\bl
The restriction $\t^{40}$ of $\t^4$
satisfies
$$
	\la \th_p, \t^{41}_p(L^4_p, L^3_p, L^2_p, L^1_p)\ra 
	= i\la \th_p, [L^4, [L^3, [L^2, L^1]]]_p\ra
$$
whenever $L^2, \1{L^1}\in H^{10}$ are in the Levi kernel
up to $\C K$-order $1$ at $p$, $L^3, L^4\in \C H$,
and $\th\in \Om^0$ is any contact form.
\el

Here we use the ``microlocal'' 
variant of the containment condition
in the Levi kernel
from
Definition~\ref{ker-1}.

\br
It is easy to see that pseudoconvexity of $M$ implies that the quartic polynomial 
$\phi^4$ in \eqref{phi24} is plurisubharmonic. 
Conversely, every plurisubharmonic $\phi^4$ appears
in a normal form of some pseudoconvex hypersurface,
e.g.\ the model hypersurface 
$$
	w+ \w = \sum_{j=1}^q |z_{2j}|^2  + \phi^4(z_4, \z_4).
$$
Furthermore, taking averages along circles,
it is easy to see that plurisubharmonicity of $\phi^4$ 
implies that of its bidegree $(2,2)$ component $\phi^{22}$.
\er

\subsection{Normal form for vector fields in the Levi kernel up to order $1$}
In any normal form as in  \eqref{phi24},
our special vector fields 
that are in the Levi kernel up to order $1$,
have particularly simple weighted expansion.
In fact, we obtain this conclusion
under a slightly more general assumption
that only requires to differentiate the Levi form
in the directions of the Levi kernel.
As customary, we shall assign weight $-a$
to coordinate vector field $\d_{z_j}$
whenever the weight of the coordinate $z_j$ is $a$.

\bp\Label{normal-vf}
Let $M$ be of the form \eqref{phi24}
and $L\in H^{10}$ be any vector field
such that 
\beq\Label{1-ker''}
	\la \th, [L^1, \1L]\ra_0
	= (L^2 \la \th, [L^1, \1L] \ra)_0
	= 0
\eeq
holds
for any $\theta\in \Om^0$,
$L^1\in H^{10}$ and 
$L^2 \in \C H$ with $L^2_0\in \C K_0$.
Then in the given coordinates and weights as 
in Proposition~\ref{4-normal},
the vector field $L$ must have 
weight at least $-1/4$ and
there exists vector $v\in \{0\}\times \C^{n-q-1}$ 
such that $L$ has a weighted expansion
\beq\Label{special-exp}
	L = \sum_j a_j ( \d_{z_{4j}} +  \phi^4_{z_{4j}} (\d_{\w} - \d_w)) +  O_w(0),
	\quad a_j\in \C,
	\quad j = 1, \ldots, n-q-1.
\eeq
In particular, $L$
cannot have any other terms of weight $-1/4$ such as
\beq\Label{low-wt}
	z_{4j}\d_{z_{2k}},
	\quad
	\z_{4j}\d_{z_{2k}}.
\eeq

Vice versa, any vector field \eqref{special-exp}
satsifies \eqref{1-ker''}.
\ep

\bpf
Since $L$ is in the Levi kernel at $0$,
its expansion cannot have any vector fields
$\d_{z^2_k}$, hence $L$ must have weight $\ge -1/4$.

To show \eqref{special-exp},
it suffices 
to prove that none of the terms \eqref{low-wt}
can occur in the expansion of $L$.
But the latter fact is a direct consquence of \eqref{1-ker''}
with 
$$
	L^1 = \d_{z^2_k}
	\quad
	L^2 = \d_{z_{4j}}
	\text{ or }
	L^2 = \d_{\z_{4j}}
	\mod \d_w.
$$

Finally, assume $L$ has expansion \eqref{special-exp}.
In particular, $L_0\in K^{10}_0$ must hold, hence
the first expression in \eqref{1-ker''} vanishes.
To show that also the second expression vanishes,
in view of Remark~\ref{kernel-span},
it suffices to assume 
that $L^1$ is either in the Levi kernel up to order $1$
or is of the form
\beq\Label{z2}
	L^1 = \d_{z_{2j}} + \z_{2j} \d_w.
\eeq
In the first case, $L^1$ also has a weighted expansion
similar to \eqref{special-exp}.
Since $\theta=i\d\rho$ has weight $\ge 1$,
a direct calculation shows that
$[L^1, \1L]$ has weight $\ge -1/4$,
and hence
$L^2 \la \th, [L^1, \1L] \ra$
 has weight $\ge -1/4 + 1 - 1/4 = 1/2$,
and therefore must vanish at $0$.

In the second case, when $L^1$ is given by \eqref{z2},
a direct calculation shows that
$[L^1, \1L]$ has weight $\ge -1/2$
and hence
$L^2 \la \th, [L^1, \1L] \ra$
 has weight $\ge -1/4 + 1 - 1/2 = 1/4$,
and therefore again must vanish at $0$.
\epf

As a byproduct,
we obtain the following
consequence:

\bc\Label{involutive-precise}
Let $L$ and $L'$ be two vector fields satisfying the assumptions
of Proposition~\ref{normal-vf}. Then
$$
	[L, L']_0, \, [\1L, L']_0 \in \C K_0.
$$
\ec

\bpf
In view of Proposition~\ref{normal-vf},
$L$, $L'$ and their conjugates commute
in their components of weight $-1/4$.
Hence their brackets must have weight $\ge -1/4$
and the statement follows.
\epf

Recall that in Definition~\ref{q-sheaf}
we introduced invariant submodule sheaves
$\6S^{10}(q)\subset H^{10}$
by requiring the condition
to be in the Levi kernel up to order $1$
to hold at every point of Levi rank $\le q$.
Then as direct application of Corollary~\ref{involutive-precise}, 
 we obtain:

\bc\Label{involutive-sheaf}
The invariant submodule sheaf 
$\6S(q) = \6S^{10}(q) \oplus \1{\6S^{10}(q) }$
satisfies the following formal integrability condition
at all points $q$ of Levi rank $\le q$:
$$
	[\6S(q), \6S(q)]_p \in \C K_p.
$$
\ec

\bpf
The statement follows by applying 
Corollary~\ref{involutive-precise}
at each point of Levi rank $\le q$.
\epf

\subsection{Symmetric extension}
Similarly to Lemma~\ref{3-symm}, we obtain a symmetric extension
for the Levi kernel restriction of $\t^4$:
\bl\Label{4-symm}
The restriction 
$$
	\t^{40}_p \colon \C K_p\times  \C K_p \times K^{10}_p\times \1{K^{10}_p} \to \C Q_p
$$ 
of the quartic tensor $\t^4_p$
 admits an unique symmetric extension
$$
	\2\t^{40}_p \colon
	\C K_p \times \C K_p\times \C K_p \times \C K_p \to \C Q_p,
$$
satisfying
\beq\Label{tensor-dif}
	\la \th_0,  \2\t^{40}_0(v^4, v^3, v^2, v^1) \ra
	= \d_{v^4} \d_{v^3} \d_{v^2} \d_{v^1} \phi^4,
\eeq
whenever $M$ is in a normal form $\rho= -2\Re w +\phi=0$
as in Proposition~\ref{4-normal}
and $\th = i\d\rho$.
In fact, \eqref{tensor-dif}
holds whenever
$\phi$
satisfies $d\phi_0=0$
and 
$ \d_{v^j} \d_{v^2} \d_{v^1} \phi^3=0$
for $j=3,4$.
\el

%
%

\section{Applications and properties of the quartic tensor}

\subsection{Relation with the D'Angelo finite type}
The quartic tensor $\t^4$ can be used to completely characterize
the finite type up to $4$ in the sense of D'Angelo~\cite{D82}
(see also ``Property P'' in \cite[Definition~5.1]{D82}):
%

\bp\Label{type-quartic}
Let $M$ be a pseudoconvex hypersurface
with nontrivial Levi kernel at $p$.
Then $M$ is of D'Angelo type $4$ at $p$
if and only if for every nonzero vector $v\in K^{10}_p$, the tensor $\t^4_p$
does not vanish when restricted to 
\beq\Label{v-res}
	(\C v + \C \1v) \times (\C v + \C \1v) \times \C v\times \C \1v.
\eeq
In fact, the latter property implies the following stronger nonvanishing condition:
$$
	\t^4_p(v,\1v, v, \1v) \ne 0.
$$
\ep

\bpf
We may assume $M$ is put into its normal form as in Proposition~\ref{4-normal}.

If the restriction of $\t^4_p$ vanishes on \eqref{v-res} for some $v\ne 0$,
we may assume $v=\d_{z_{31}}$, where $z_3=(z_{31},\ldots, z_{3,n-r})$.
Then it follows from the normal form that the line $\C v$ has order of contact with $M$
higher than $4$, hence the D'Angelo type at $p$ is also higher than $4$.

On the other hand, suppose the restriction of $\t^4_p$ to \eqref{v-res} does not vanish for any 
$v\ne 0$. Assume by contradiction, 
there exists a nontrivial holomorphic curve 
$$
	\g\colon (\C,0) \to (\C^{n+1},0),
	\quad
	\g(t) = \sum_{k\ge k_0} a_k t^k,
	\quad 
	a_{k_0}\ne 0,
$$
whose contact order with $M$ at $0$ is higher than $4$.
Recall that the contact order is given by
$$
	\frac{\nu(\rho\circ\g)}{\nu(\g)},
$$
where $\rho$ is any defining function of $M$ and $\nu$ is the vanishing order at $0$,
in particular, $\nu(\g) = k_0\ge 1$.
Taking $\rho := -2\Re w +\phi$, we must have $a_{k_0}\in \{0\} \times \C^n$,
otherwise the contact order would be $1$.
Similarly, expanding $\rho\circ \g$, it follows by induction that 
$$
	a_{l}\in \{0\} \times \C^{n},
	\quad
	l< 4k_0,
$$
and
$$
	a_{k}\in \{0\} \times \{0\} \times \C^{n-r},
	\quad
	k< 2k_0,
$$
for otherwise the contact order would be less than $4$.
Finally collecting terms of order $4k_0$
and using our assumption that the contact order is greater than $4$,
we obtain
\beq\Label{4k}
 	\phi^2(a_{2k_0} t^{2k_0}, \1a_{2k_0} \1t^{2k_0} ) + \phi^4(a_{k_0} t^{k_0}, \1a_{k_0}\1t^{k_0}) =0.
\eeq
In particular, it follows that 
$$
	\phi^4(a_{k_0} \xi, \1a_{k_0}\1\xi) = c\xi^2\bar\xi^2.
$$
Since $\phi^4$ is plurisubharmonic, we must have $c\ge 0$.
Hence both terms in \eqref{4k} are nonnegative,
and therefore must vanish.
In particular, 
$
	c a_{k_0}^2 \1a_{k_0}^2 = 0,
$
implying 
$$
	\phi^4(a_{k_0} \xi, \1a_{k_0}\1\xi) = 0,
$$
which is in contradiction with our nonvanishing assumption on $\t^4_p$.
Hence the D'Angelo type is $4$ completing the proof of the converse direction.

Finally, the last statement follows from the plurisubharmonicity of $\phi^4$
in any normal form.
\epf

The last conclusion suggests a connection
with the so-called regular type.
Recall that the {\em regular type}
of $M$ at $p$
is the maximum (possibly infinite) of the vanishing order 
$\nu(\rho\circ\g)$,
where $\rho$ is a defining function of $M$ in a neighborhood of $p$
and $\g\colon (\C,0)\to (\C^n,p)$ a germ
of a regular complex-analytic curve,
i.e.\ satisfying $\g'(0)\ne 0$.
As a consequence of Proposition~\ref{type-quartic},
we obtain the following characterization:

\bc
Let $M$ be a pseudoconvex hypersurface
with nontrivial Levi kernel at $p$.
Then the following are equivalent:
\begin{enumerate}
\item $M$ is of D'Angelo type $4$ at $p$;
\item $M$ is of regular type $4$ at $p$;
\item the quartic tensor satisfies $\t^4_p(v,\1v, v, \1v) \ne 0$
whenever $v$ is a nonzero vector in the Levi kernel at $p$.
\end{enumerate}
\ec

\bpf
The equivalence of (1) and (3) is contained in Proposition~\ref{type-quartic}. The equivalence of (2) and (3) is obtained by repeating the proof of the proposition for a regular curve. 
\epf

The pseudoconvexity assumption in Proposition~\ref{type-quartic} cannot be dropped:

\be
Let $M\subset \C^3_{w, z_1, z_2}$ be given by
$$
	2\Re w = |z_1|^2 - |z_2|^4.
$$
Then $M$ contains the image of the curve $t\mapsto (0, t^2, t)$ 
and is hence of infinite type at $0$.
On the other hand, $M$ is in the normal form \eqref{phi24}
and hence $\t^4_0(v, \1v, v, \1v)\ne 0$  for any $v\ne0\in K^{10}_0$.
\ee

%

\subsection{Uniformity of the quartic tensor}
The sheaves  $\6S^{10}(q)$
introduced in Definition~\ref{q-sheaf}
 can be used to 
obtain a uniform behavior of 
$\t^4_p$ as $p$ varies
over the set of nearby points
of bounded Levi rank.
In fact, as direct consequence
from the definition and Corollary~\ref{span-ker}, 
we obtain
that  $\t^4_p$
can be calculated using local sections of $\6S^{10}(q)$:
\bc\Label{uni-t4}
For every 
vector fields $L^4, L^3\in \C T$
and $L^2, \1{L^1}\in \6S^{10}(q)$
defined in an open set $U\subset M$,
the identity 
\eqref{1-levi'} holds
simultaneously
for all points $p\in U$
of Levi rank $q$.
\ec

\br
In the context of Corollary~\ref{uni-t4},
it is essential to require the vector fields $L^2, L^1$
to be contained in Levi kernels up order $1$ 
(rather than merely contained in Levi kernels).
In fact, for a higher order perturbation of Examples~\ref{ex-4} where 
$0$ is the only Levi-degenerate point, 
choosing vector fields $L^j$ as higher order perturbations of
the vector field $L$ in the example or its conjugate
would violate \eqref{1-levi'}.
\er

It is important to note that the conclusion of Corollary~\ref{uni-t4}
may not hold for points $p\in U$ of Levi rank $>q$
when $L^2_p, \1{L^1_p}\in K^{10}_p$.
In fact, $\t^4_p$ may not even be continuous 
e.g. may vanish for $p$ of higher Levi rank even when $\t^4_{p_0}$ does not vanish 
on any line for $p_0$ of Levi rank $q$.
This is illustrated by D'Angelo's celebrated example 
where the finite type is not upper-semicontinuous \cite{D80, D82},
see also Example~4 and its continuation
on pages 135--136 in \cite{D93}:

\be[D'Angelo]\Label{da-ex}
Let $M\subset\C^3_{w, z_1, z_2}$ be given by
$$
	2\Re w = |z_1^2 - w z_2|^2 + |z_2|^4.
$$
Then $M$ is of Levi rank $0$ and finite type $4$ at $0$ and hence $\t^4_0$ does not vanish on the lines 
products \eqref{v-res} in view of Proposition~\ref{type-quartic}.
In fact, $M$ is in its normal form as in Proposition~\ref{4-normal} with
$\phi^4 = |z_1|^4 + |z_2|^4$, and hence 
$$
	\t^4_0(v, \1v, v, \1v) = 4(|v_1|^4 + |v_2|^4),
	\quad v\in K^{10}_0 \cong \{0\} \times C^2_{z_1, z_2}.
$$

On the other hand, at every $p=(it, 0, 0)$ on the imaginary axis with $t\ne 0$,
the Levi rank is $1$, and $M$ can be locally transformed into a normal form 
\eqref{phi24} with vanishing $\phi^4$ implying $\t^4_p(v, \1v, v, \1v)=0$
for any $v\in K^{10}_p$. Thus $\t^4_p(v, \1v, v, \1v)$ cannot be continuous
for any $v=v(p)$ converging to any $v(0)\ne 0$
as $p\to 0$.

Of course, this phenomenon is closely related to the lack of
 upper-semicontinuity of the type as demonstrated by D'Angelo.
The additional importance of this example and its generalisation
in \cite[Example 5.16]{D82} and 
\cite[Example 4]{D93} is the occurence of the 
``worst possible lack of semicontinuity'' of the type
for pseudoconvex hypersurfaces,
demonstrating sharpenss of the D'Angelo's bound
controlling the type in a neighborhood of a point $p\in M$
in terms of the type at $p$, see \cite[Theorem 5.5]{D82}.
\ee

\subsection{Kernels of quartic tensors}
For any homogenous polynomial, consider the following notion of
holomorphic kernel:
\bd
	The {\em holomorphic kernel of a homogeneous polynomial} 
	$P(z,\bar z)$, $z\in \C^n$, 
	is defined to be the subspace of all $(1,0)$ vectors $v$ such that
	\beq\Label{kernel-def}
		\d_v P(z,\z) \equiv \d_{\1v}P(z,\z) \equiv 0.
	\eeq
	Equivalently, the holomorphic kernel is the space of all $v$
	such that both $v$ and $\1v$ belong to the kernel 
	of the polarization of $p$.
\ed
It is straightforward to see the following simple characterization of the kernel:
\bl\Label{ker-coor}
	The holomorphic kernel of $p$ is the maximal subspace $V$ such that,
	there exists a linear change of coordinates such that
	$$V= \oplus_{j=1}^l (\C \d_{z_j} \oplus \C \d_{\z_j})$$
	and $P(z,\z)$ is independent of the variables $z_1,\ldots, z_l$ and their conjugates.
\el

\bd\Label{rank}
	The {\em rank} of $P$ is $n-d$, where $d$ is the dimension of the holomorphic kernel.
\ed

Also separating bihomogeneous components in \eqref{kernel-def}, we obtain:
\bl
	Let 
	$$
		P(z,\z)=\sum P_{kl}(z,\z)
	$$ 
	be a decomposition into components $P_{kl}$ of bidegree $(k,l)$ in $(z,\z)$.
	Then the holomorphic kernel of $P$ equals the intersection of kernels of $P_{kl}$ for all $k,l$.
\el

Next we compare the holomorphic kernel of the polynomial $\phi^4$
in the normal form given by Proposition~\ref{4-normal}
and the restriction 
$$
	\t^{40}_p \colon \C K_p \times \C K_p \times K^{10}_p \times \1{K^{10}_p} \to \C Q_p.
$$
of the quartic tensor $\t^4_p$
to the Levi kernel in each component.

\bd
The {\em holomorphic kernel} of $\t^{40}_p$ is 
$V\cap \1V$, where
$$
	V = \ker \t^{40}_p = \{ v\in \C K_p : 
		\t^{40}_p(v, v^3, v^2, v^1) =0
		\text{ for all }
		v^3, v^2, v^1\}.
$$
\ed

First of all, remark that 
without pseudoconvexity assumption,
the holomorphic kernel of $\t^{40}_p$
may get larger than that of $\phi^4$:
\be\Label{diff-kernels}
Let $M\subset \C^3_{w, z_1, z_2}$ be given by
$$
	2\Re w = \phi^4(z,\z):= 2\Re (z_1^3 \z_2).
$$
Then the arguments in the proof of Proposition~\ref{4-normal}
can be used to show that \eqref{phi-diff} still holds,
implying that $\d_{z_2}$ and $\d_{\z_2}$
are in the kernel of $\t^{40}_0$
in the $1$st and $2$nd arguments
but not in the $3$rd one.
\ee

On the other hand, in presence of pseudoconvexity,
both kernels must coincide as the following lemma shows.
As a matter of convention, for a multilinear function 
$f(v^1, \ldots, v^m)$, we call its kernel in the $k$th argument
the space of all $v^k$ such that $f(v^1, \ldots, v^m)=0$ holds
for all $v^j$ with $j\ne k$.

\bl\Label{ker-rel}
	Let $M$ be in its normal form given by Proposition~\ref{4-normal},
	and
	assume that $M$ is pseudoconvex.
	Then both holomorphic kernels of $\t^{40}$ in the $1$st and $2$nd arguments coincide with holomorphic kernel 
	$V$ of $\phi^4$.
	Furthermore,
	the kernels of $\t^{40}$ in the $3$rd and $4$th arguments
	 coincide respectively with $V$ and $\1V$.
\el

\bpf
As direct consequence of \eqref{phi-diff} we obtain
that the holomorphic kernel of $\phi^4$ is contained
in the kernel of $\t^{40}_p$ in each argument.

Vice versa, let $v$ be $(1,0)$ vector in the holomoprhic 
kernel of $\t^{40}_p$ (in the $1$st argument).
We write $\xi=z_4$ for brevity.
After a linear change of coordinates we may assume 
$v=\d_{\xi_1}$, where $\xi_1$ is the first component of $\xi$
in the notation of Proposition~\ref{4-normal}.
Then it follows from \eqref{phi-diff} that
$\d_{\xi_1} \phi^4$ is harmonic.
Since $\phi^4$ has no harmonic terms, it must have the form
$$
	\phi^4 = 2\Re (\bar\xi_1 h(\xi)) + R,
$$
where $h$ is holomorphic and $R$ is independent of $\xi_1$.
Now since $M$ is pseudoconvex, $\phi^4$ is plurisubharmonic,
in particular,
\beq\Label{t}
	(\d_{\xi_1} + t \d_{\xi_j}) (\d_{\bar\xi_1} + t \d_{\bar\xi_j}) 
	\phi^4 \ge 0
\eeq
holds for all $t\in \R$.
Then for $t=0$, we obtain $\d_{\xi_1}h\ge 0$.
Since $h$ is holomorphic, we must have $\d_{\xi_1}h \equiv 0$.
Hence the linear part of \eqref{t} must be $\ge 0$
and therefore equal to $0$, since $t$ is any real number.
But this means $h\equiv 0$, and hence $v=\d_{\xi_1}$ is in the 
holomorphic kernel of $\phi^4$ as claimed.

The claimed statements for kernels in other arguments of $\t^{40}_p$
are obtained by 
repeating the same proof. 
\epf

In view of Lemma~\ref{ker-rel},
we simply refer to the {\em holomorphc kernel of $\t^{40}$}
for its kernel in the $1$st (and, equivalently, in the $2$nd) argument,
Also the {\em rank of $\t^{40}$} is 
$\dim K^{10}_p - d$,
where $d$ is the dimension of its holomorphic kernel,
which coincides with the rank of $\phi^4$
in the sense of Definition~\ref{rank}.

\section{Relation with Catlin's multitype}
\Label{m-type}

We first recall the definition of the multitype due to Catlin~\cite{C84proc}.
Consider 
arbitrary ordered weights
\beq\Label{tuple}
	\L = (\l_1, \ldots, \l_{n}), 
	\quad
	1\le \l_1 \le \ldots \le \l_{n} \le \infty,
\eeq
and for a multiindex $\a=(\a_1,\ldots, \a_n)$,
define its weight by
$$
	\|\a\|_\L := \l_1^{-1} \a_1 +\ldots + \l_n^{-1} \a_n.
$$

\bd\Label{adm}
A weight $\L$ is called {\em admissible} if
for each $k=1,\ldots,n$, either $\l_k=+\infty$ or there exists 
a $k$-tuple of nonnegative integers $a=(a_1,\ldots, a_k)$
satisfying $a_k>0$ and $\|a\|_\L = 1$.
\ed

Next, for a smooth function $\rho(z,\z)$ defined in a neighborhood of $0$, write
\beq\Label{rho-wt}
	\rho = O_{\L}(1)
\eeq
whenever all nonzero monomials $\rho_{z^\a\z^\b}z^\a\z^\b$ in the Taylor expansion of $\rho$
at $0$ satisfy $\|\a+\b\|_\L \ge 1$,
i.e.\ are of weight $\ge 1$,
or, equivalently, whenever the estimate
$$
	|\rho(z,\z)| \le C (|z_1|^{\l_1} + \ldots + |z_n|^{\l_n})
$$
holds in a neighborhood of $0$ for some $C> 0$.
Finally, we regard $(\mu_1, \ldots, \mu_n)$
as lexicographically smaller
than $(\l_1, \ldots, \l_{n})$
whenever $\mu_j <\l_j$
holds for the smallest $j$ such that $\mu_j\ne \l_j$.


\bd[\cite{C84proc}]\Label{mt}
The {\em multitype} of a smooth real hypersurface $M\subset \C^{n}$ at 
a point $p\in M$
is the lexicographic supremum of the set of 
all admissible weights $\L$
such that there exist local 
holomorphic coordinates in a neighborhood of $p$
and vanishing at $p$,
where $M$ is given by
a defining function $\rho$  
satisfying 
$\rho = O_\L(1)$. 
\ed

Note that since a defining function $\rho$ is unique
up to a nonvanishing smooth factor,
the property $\rho = O_\L(1)$
is independent of the choice of $\rho$
(but, of course, in general, it does depend on the coordinates).

Catlin's multitype is an essential ingredient in his proof
of the Property (P) \cite{C84ann}, based on \cite{C84proc},
required for the proof of global regularity of the $\bar\d$-Neumann problem,
as well as for its quantitative analogue
in \cite{C87} required for the proof of subelliptic estimates.
However, in practice, the multitype is difficult to compute,
due to the nature of its definition requiring
taking lexicographic supremum over arbitrary holomorphic coordinate charts.

In \cite{Ko10}, Kolar gave a remarkable general algorithm
for computing the multitype of a given hypersurface,
which involves taking certain consecutive coordinate normalisations.
However, the number of steps involved there equals the dimension,
and each step depends on the previous coordinate choice.
%
%
%
%
In our case, we can use the quartic tensor
$\t_4$ to have an algebraically invariant way of calculating part 
of the multitype, that works regardless of the dimension $n$
and independent of coordinates.
In case $\t_4$ has trivial kernel, this approach gives the complete multitype.

\bp\Label{multi-quartic}
Let $M\subset \C^{n}$ be a pseudoconvex hypersurface,
and $p\in M$ a point with Levi form of rank $q_2$ 
and the restricted quartic tensor $\t^{40}$ of rank $q_4$.
Then the multitype $\L = (\l_1, \ldots, \l_{n})$ of $M$ at $p$
satisfies
\beq\Label{wt-eq}
	\l_1 = 1, 
	\quad
	\l_2 = \ldots \l_{q_2+1} =2, 
	\quad
	\l_{q_2+2} = \ldots = \l_{q_2 + q_4 +1} = 4,
\eeq
and 	
\beq\Label{wt-ineq}
	\l_k > 4, \quad k> q_2 + q_4 +1.
\eeq

In particular, if $\t^{40}$ has only trivial kernel,
the multitype is
$(1, 2,\ldots, 2, 4, \ldots, 4)$,
where the number of $2$'s equals the Levi rank.
\ep

\bpf
By Lemma~\ref{ker-coor}, 
in addition to the normal form in Proposition~\ref{4-normal},
we can make $\phi^4$ independent of the last $d$ coordinates,
where $d$ is the dimension of the $(1,0)$ kernel of $\t_4$.
This shows that it is possible to achieve 
\eqref{rho-wt}
with weights
satisfying both \eqref{wt-eq} and \eqref{wt-ineq}.

The actual multitype may only be lexicographically higher,
in particular, \eqref{wt-ineq} is already satisfied.
Assume by contradiction that we have another choice of coordinates 
with higher weights failing one of the equalities in \eqref{wt-eq}.
However, we must obviously have $\l_1=1$
and the Levi form invariance forces the next $q_2$ weights to be equal $2$.
Therefore we must have some $\l_k > 4$ for $k\le q_2+q_4+1$.
In those coordinates, we would have the same normal form as in Proposition~\ref{4-normal}
with $\phi^4$ being independent of $z_j$ at least for $j\ge q_2 + q_4 +1$.
That, however, would mean that the rank of $\t_4$ is less than $q_4$, which is a contradiction.
Hence the multitype must satisfy all of \eqref{wt-eq} as claimed,
where the admissibility condition from Definition~\ref{adm} is clearly satisfied.
\epf

\section{Ideal sheaves for the Levi rank level sets}\Label{ideal-par}

We use the vector field submodule 
sheaves $\6S^{10}(q)$ in Definition~\ref{q-sheaf}
to define invariant ideal sheaves of smooth functions 
$\6I(q)$ (as in Theorem~\ref{main}, part (5)):

%
%
%
%

\bd\Label{ideal}
Let $M\subset\cn$ be a pseudoconvex hypersurface.
For every $q$,
define $\6I(q)$ to be the ideal sheaf 
generated by all (smooth complex) functions $g$,$f$ of the form
$$
	g = \la \th, [L^2, L^1] \ra,
	\quad
	f = L^3 \la \th, [L^2, L^1] \ra,
$$
where $\th\in\Om^0$ is a contact form, 
$L^3\in \C T$ arbitrary complex vector field,
and $L^2, \1{L^1} \in \6S^{10}(q)$
arbitrary sections.
\ed

\br
The same ideal sheaf $\6I(q)$ is generated
by the functions $g,f$ as in Definition~\ref{ideal},
where the vector fields $L^j$ can be chosen
from a fixed finite set of generators
of $\C T$ and $\6S^{10}(q)$.
That is due to the linearity of $g$
with respect to smooth functions,
and the linearity of $f$ modulo
the ideal generated by all functions $g$.
\er

\be
Example~\ref{ex-module-0}
shows that for $q=0$,
the sheaf $\6S^{10}(0)$
contains all germs of all $(1,0)$ vector fields.
Then the ideal sheaf $\6I(0)$
is generated by all Levi form entries
and their first order derivatives.

On the opposite end,
Example~\ref{ex-module-0}
shows that for $q=n-1$,
the sheaf $\6S^{10}(n-1)$
consists of all germs of $(1,0)$ vector fields
that are everywhere contained in the Levi kernel.
Such sheaf is always trivial when
the hypersurface $M$ is generically Levi-nondegenerate,
which is the case e.g.\ whenever $M$ is of finite type.
In the latter case, the ideal sheaf $\6I(n-1)$ is also trivial
(identically zero),
which corresponds to the simple fact
that the set of points of the maximal Levi rank $n-1$
is never contained in a proper submanifold.
\ee

\be
Let $M$ be as in Example~\ref{z12}.
Then away from the subset $M^1 = \{z_1z_2=0\}\subset M$
of the points of Levi rank $\le 1$,
the ideal sheaf $\6I(1)$ is generated by the Levi form entries
$|z_1|^2, |z_2|^2$, 
that generate 
all smooth germs of functions there.
This is expected as that set consists
of the Levi-nondegenerate points.

On the other hand, along the Levi-degeneracy set $M^1$,
the Levi form and its derivatives in Definition~\ref{ideal}
need to be computed along vector fields 
$L^2, \1{L^1}$
from the submodule $\6S^{10}(1)$.
In view of Example~\ref{z12},
the Levi form entries
$g = \la \th, [L^2, L^1] \ra$
always vanish of order at least $2$ along $M^1$,
whereas their derivatives
$f = L^3 \la \th, [L^2, L^1] \ra$
generate the maximal ideal of $M^1$ away from
$z_1=z_2=0$.
\ee

As a direct consequence of Corollary~\ref{span-ker} and 
Lemma~\ref{psc-vanish}, we obtain
a general way of constructing submanifolds
containing level sets of the Levi rank:

\bc\Label{iq-van}
Let $M$ be a pseudoconvex hypersurface.
Then every local section in $\6I(q)$
vanishes at all points of Levi rank $q$.
In particular, for any collection $f^1, \ldots, f^m$
of real functions from the real part $\Re \6I(q)$
defined in an open set $U\subset M$
satisfying 
\beq\Label{wedge}
	df^1 \wedge \ldots \wedge df^m\ne 0,
\eeq
the submanifold
$$
	S = \{ f^1 = \ldots = f^m =0 \}
$$
contains the set of all points of Levi rank $q$ in $U$.
In fact, the set $S$ still has the same property
without assuming \eqref{wedge}.
\ec

\br
Note that due to our definition of $\6I(q)$,
any complex multiple of a local section is again a local section.
Consequently, it suffices to take only sections in $\Re \6I(q)$
to define the same set.
\er

We next apply the quartic tensor to describe the differentials
of sections in $\6I(q)$.

\bd
For an ideal sheaf $\6I$ define it {\em kernel} at $p$
$$
	\ker_p \6I \subset \C T_p,
$$
to be the intersection of kernels of all differentials $df_p$,
where $f$ is any local section of $\6I$ in a neighborhood of $p$.
\ed

Then Corollary~\ref{iq-van} implies:

\bl
Let $p\in M$ be a point of Levi rank $q$.
Then the kernel of the $q$th sublevel ideal $\6I(q)$ at $p$
coincides with the kernel of the quartic tensor $\t^4_p$.
\el

We can now summarize this paragraph's results as follows:

\bp\Label{main0}
Let $M\subset\cn$ be a pseudoconvex real hypersurface,
and $p\in M$ a point of Levi rank $q$.
Then in a neighborhood of $p$, 
the set 
 of all points of the same Levi rank $q$
is contained in a real submanifold $S\subset M$
through $p$ such that
$$
	T_p S = \ker \t^4_p,
$$
and $S$ is given by the vanishing 
of local sections
$$
	f^1, \ldots, f^m \in \6I(q),
	\quad
	df^1\wedge \ldots \wedge df^m \ne 0.
$$
%
In particular, when $M$ of finite type $4$ at $p$,
the intersection of $T_pS$
with the Levi kernel at $p$ is totally real.
\ep

%
%


\section{Relation with Catlin's boundary systems}\Label{bs}

\subsection{Maximal Levi-nondegenerate subbundles}


Recall that {\em Catlin's boundary system} construction
for a hypersurface $M$ at a point $p\in M$
begins with a maximal collection of $(1,0)$ vector fields $L_2, \ldots, L_{q+1}$ tangent to $M$
such that the Levi form matrix 
$$
	(\la \th,  [L_j, \1L_k] \ra)_{2\le j,k \le q}
$$
is nonsingular. In particular, $q$ must be equal to the Levi rank at $p$.

Invariantly, consider any {\em maximal Levi-nondegenerate subbundle} through $p$,
i.e.\ any smooth subbundle  $V^{10} \subset H^{10}$ where the restriction of the Levi form is nondegenerate. 
Then obviously any such $V^{10}$ appears as the span
of the first $q$ vector fields in Catlin's boundary system,
and vice versa, every such span is a maximal Levi-nondegenerate subbundle.
%

Next Catlin considers the Levi-orthogonal subbundle 
$$
	S^{10} := (V^{10})^\perp \subset H^{10}
$$ 
($T^{10}_{q+2}$ in Catlin's notation).
In particular, the subbundle {\em $S^{10}$ contains all Levi kernels $K^{10}_x$ at all points $x\in M$ near $p$},
even when $\dim K^{10}_x$ depends on $x$.
That makes the fiber $(V^{10}_x)^\perp$ unique
whenever the Levi rank at $x$ is the same as $p$,
even when $V^{10}_x$ itself may not be unique.
On the other hand, at points $x$ of higher Levi rank,
$(V^{10}_x)^\perp$ clearly depends on the choice of $V^{10}_x$.

The rest of Catlin's boundary system construction 
only depends on the subbundle $S^{10}$ 
rather than on $V^{10}$ and its chosen basis.

\subsection{Levi kernel inclusion of higher order}
As mentioned before, $S^{10}$ contains the Levi kernel at every point.
On the other hand, if $M$ is pseudoconvex,
we have shown in Lemma~\ref{1-ker-def} 
that $S^{10}$ is itself 
{\em contained in the Levi kernel up to order $1$ at $p$}
as defined in Definition~\ref{ker-1}.
That permits to use arbitrary sections of $S^{10}$
in the calculation of the quartic tensor $\t^4$:

\bc\Label{s10-tensor}
Let $M$ be a pseudoconvex hypersurface,
$V^{10}\subset H^{10}$
a maximal Levi-nondegenerate subbundle at $p\in M$,
and $S^{10}$ the Levi-orthogonal complement of $V^{10}$.
Then the quartic tensor $\t^4_p$ defined by Lemma~\ref{d2}
satisfies 
$$
	\la \th_p, \t^4_p(L^4_p, L^3_p, L^2_p, L^1_p)\ra 
	= i(L^4 L^3 \la \th, [L^2, L^1]\ra)_p
$$
for any $L^4, L^3 \in \C T$, $L^2, \1{L^1}\in S^{10}$ and $\th\in \Om^0$.
\ec




%
%

\subsection{Relation with the rest of Catlin's boundary system construction}
\Label{bs-2}
The remaining part of Catlin's construction is based on the higher order Levi form derivatives
\beq\Label{L-th}
	\6L \th := L^m \ldots L^3 \la \th, [L^2, L^1] \ra,
	\quad
	\6L = (L^m, \ldots, L^1),
\eeq
where $\th = \d r$ and $r$ is a defining function of $M$.
Then a boundary system 
\beq\Label{bd-sys}
	\6B = \{r_1, r_{q_2+2}, \ldots, r_\nu; L_2, \ldots, L_\nu  \},
	\quad q_2+2 \le \nu\le n,
\eeq
is constructed
together with associated weights
$$
	\a_1 =1 <  \a_2 =\ldots = \a_{q}=2 < \a_{q+1}\le \ldots \le \a_{\nu}\le \infty,
$$
where $r_1 = r$ is the given defining function,
$L_j$ and $r_j$ are respectively smooth $(1,0)$ vector fields
and smooth real functions in a neighborhood of $p$.
The construction proceeds by induction as follows.
Assuming a boundary system is constructed for given $\nu$,
define the next subbundle
$$
	T^{10}_{\nu+1} := \{ L\in T^{10}_{q_2+2} : \d r_{q_2+2}(L) =\ldots = \d r_{\nu}(L) =0 \}.
$$
Then count all previous $L_j$ and their conjugates with weight $\a_j$,
and consider a new vector field $L_{\nu+1}\in T^{10}_{\nu+1}$ and its conjugate,
whose weight $\a=\a_{\nu+1}$ is to be determined.
Now look for all lists 
$\6L = (L^m, \ldots, L^1)$ with each $L^k\in \{L_{q+2}, \ldots, L_{\nu+1}\}$,
which are of total weight $1$ and {\em ordered}, 
i.e.\ $L_j, \1L_j$ preceed $L_k, \1L_k$ whenever $j>k$,
such that 
\beq\Label{non-vanish}
	(\6L\d \rho)_p\ne 0.
\eeq
The list must contain the new vector field $L_{\nu+1}$ or its conjugate,
and the new weight $\a_{\nu+1}$ is chosen to be minimal possible
with that property.
Finally set either
$$
	r_{\nu+1}:= \Re L^{m-1}\ldots L^3 \la \th, [L^2, L^1] \ra
	\text{ or }
	r_{\nu+1}:= \Im L^{m-1}\ldots L^3 \la \th, [L^2, L^1] \ra	
$$
such that 
$$
	(L_{\nu+1} r_{\nu+1})_p\ne 0,
$$
which is always possible in view of \eqref{non-vanish}, 
since the first vector field in the list,
 $L^m$ is either $L_{\nu+1}$ or its conjugate.
Restating
Lemma~\ref{1-ker-def} and
Corollariy~\ref{s10-tensor}, 
we have:

\bc
Let $M$ be pseudoconvex hypersurface 
with Levi form of rank $q$ at $p$.
Fix a boundary system $\{L_2, \ldots, L_{q+1}\}$ at $p$.
Then $S^{10}= T^{10}_{q+2} = V^\perp$ for 
$V:=\span\{L_2, \ldots, L_{q+1}\}$.
Further, for any vector fields 
$L^4, L^3\in S^{10} + \1{S^{10}}$, $L^2\in S^{10}$, $L^1\in \1{S^{10}} $,
we have
$$
	L^3 \la \th, [L^2, L^1] \ra_p = 0,
$$
$$
	L^4 L^3 \la \th, [L^2, L^1] \ra_p = \t^{40}_p(L^4_p, L^3_p, L^2_p, L^1_p).
$$
In other words, for lists $\6L$ of length $3$, 
the derivative $(\6L\th)_p$ vanishes,
whereas for lists of length $4$, 
it only depends on the vector field values at $p$
and is given by the restricted quartic tensor $\t_4$
(regardless of the choice of the boundary system).
\ec

Thus via the quartic tensor restriction $\t^{40}_p$,
the nonvanishing condition in \eqref{non-vanish}
is reduced to a purely algebraic property
only depending on the vector fields' values at $p$.



%


%
%
%
%
%

\end{document}